\documentclass[11pt]{article}
\usepackage{amssymb,amscd,latexsym}   

\newcommand{\rar}{\rightarrow}
\newcommand{\lar}{\longrightarrow}
\newcommand{\surjects}{\twoheadrightarrow}

\newtheorem{Theorem}{Theorem}[section]
\newtheorem{Lemma}[Theorem]{Lemma}
\newtheorem{Corollary}[Theorem]{Corollary}
\newtheorem{Proposition}[Theorem]{Proposition}
\newtheorem{Remark}[Theorem]{Remark}

\newtheorem{Example}[Theorem]{Example}

\newtheorem{Definition}[Theorem]{Definition}
\newtheorem{Question}[Theorem]{Question}

\def\demo{\noindent{\bf Proof. }}
\def\text{\mbox}
\def\phi{\varphi}

\newcommand{\rk}{\mbox{\rm rank}}

\newcommand{\hht}{\mbox{\rm height }}

\def\pp{\mathbb P}

\def\xx{{\bf x}}
\def\TT{{\bf T}}

\def\yy{\bf y}
\def\zz{\bf z}
\def\ff{{\bf f}}

\def\fd{{\cal D}}
\def\fp{{\cal P}}
\def\fz{{\cal Z}}

\newcommand{\lcm}{\mbox{\rm lcm }}

\def\zw{\mbox{\sc z}_{\mathfrak w}}
\def\zwt{\widetilde{\mbox{\sc z}_{\mathfrak w}}}
\def\tw{\mbox{\sc t}_{\mathfrak w}}

\def\bg{{\bf g}}
\def\hh{{\bf h}}

\def\pp{{\mathbb P}}

\def\QED{\hfill$\Box$}
\def\qed{\QED}

\begin{document}

\begin{center}

 {\bf {\LARGE Polar syzygies in characteristic zero: \\[5pt]
  the monomial case}}

\vspace{0.3in}

\footnotetext{Partially supported by {\it Ministerio de Educaci\'on
y Ciencia -- Espa\~na} (MTM2007-61444)}

{\large\sc Isabel Bermejo}\footnote[1]{The first two authors thank
the {\it Universidade Federal de Pernambuco} for hospitality and
partial support during the preparation of this work. }, {\large\sc
Philippe Gimenez}\footnotemark[1] and {\large\sc Aron
Simis}\footnote[2]{This author thanks C.I.M.A.C. for providing
support, and the {\it Universidad de La Laguna} for hospitality.}

\end{center}

\begin{abstract}
Given a set of forms ${\ff}=\{f_1,\ldots,f_m\}\subset
R=k[x_1,\ldots,x_n]$, where $k$ is a field of characteristic zero,
we focus on the first syzygy module ${\cal Z}$ of the transposed
Jacobian module ${\cal D}({\bf f})$, whose elements are called {\it
differential syzygies\/} of ${\bf f}$. There is a distinct submodule
${\cal P}\subset {\cal Z}$ coming from the polynomial relations of
${\bf f}$ through its transposed Jacobian matrix, the elements of
which are called {\it polar syzygies\/} of ${\bf f}$. We say that
${\bf f}$ is {\it polarizable\/} if equality ${\cal P}= {\cal Z}$
holds. This paper is concerned with the situation where ${\bf f}$
are monomials of degree $2$, in which case one can naturally
associate to them a graph ${\cal G}({\bf f})$ with loops and
translate the problem into a combinatorial one. A main result is a
complete combinatorial characterization of polarizability in terms
of special configurations in this graph. As a consequence, we show
that polarizability implies normality of the subalgebra $k[{\bf
f}]\subset R$ and that the converse holds provided the graph ${\cal
G}({\bf f})$ is free of certain degenerate configurations. One main
combinatorial class of polarizability is the class of polymatroidal
sets. We also prove that if the edge graph of ${\cal G}({\bf f})$
has diameter at most $2$ then ${\bf f}$ is polarizable. We establish
a curious connection with birationality of rational maps defined by
monomial quadrics.
\end{abstract}

\section{Introduction}

Let $k$ be a field of characteristic zero. Given a set of forms of
the same degree, ${\ff}=\{f_1,\ldots,f_m\}\subset
R=k[x_1,\ldots,x_n]$, one can consider both the ideal
$I=(\ff)\subset R$ and the $k$-subalgebra
$A=k[\ff]=k[f_1,\ldots,f_m]\subset R$. Looking at the intertwining
properties of the subalgebra $A$ and the ideal $I$ was of course
Hilbert's original idea to understand the finite generation of
certain rings of invariants. As such it became natural to look at
the syzygies of the polynomial relations of $I$. About 25 years
before Hilbert's wrap-up of these questions, P. Gordan and M.
Noether in their celebrated work \cite{GN} about the Hesse problem
had this approach sort of turned around by looking instead at an
individual polynomial relation $F\in k[\TT]=k[T_1,\ldots,T_m]$ of
$\ff$ in the special case where $n=m$ and $f_1,\ldots,f_m$ were the
partial derivatives of a homogeneous polynomial $f\in R$. They posed
(and solved) the question of finding all polynomial solutions
$\Phi(\xx)\in R$ of the partial differential equation
\begin{equation}\label{diffeq}
\sum_{j=1}^m \frac{\partial \Phi}{\partial x_j}\, F_{T_{j}}(\ff)=0,
\end{equation}
where a subscripted variable indicates partial derivative with
respect to this variable. In other words, among all syzygies of the
ideal $(F_{T_1}(\ff),\ldots,F_{T_m}(\ff))$ they were looking for the
polynomially integrable ones! Particular solutions are of course
 the very partial derivatives $\Phi_i=f_i$, one for each $i=1,\ldots,m$ -- a
consequence of the rule of derivatives for composite functions.

\smallskip

Now,  one can think about the relations
$$\sum_{j=1}^m \frac{\partial^2 f}
{\partial x_i\partial x_j}\,\frac{\partial F}{\partial T_j}(\ff),\;
i=1,\ldots,m,
$$
for each polynomial relation $F\in k[\TT]$ of $\ff$, as syzygies
of the Hessian matrix of the form $f$. Going back to the more
general setting where ${\ff}=\{f_1,\ldots,f_m\}$ is a set of $m$
forms of the same degree in $R=k[x_1,\ldots,x_n]$, one could ask for the syzygies
of the transposed Jacobian matrix
of ${\ff}$. This was the original goal in \cite{jac} where the
syzygies corresponding to the relations
$$\sum_{j=1}^m \frac{\partial f_j}
{\partial x_i}\,\frac{\partial F}{\partial T_j}(\ff),\;
i=1,\ldots,m,
$$
one for each polynomial relation $F\in k[\TT]$ of $\ff$, have been
dubbed {\it polar syzygies\/} and it was shown that in a certain
special context the whole module of syzygies of the transposed
Jacobian matrix of ${\ff}$ is generated by the polar syzygies.

\smallskip

The motivation for the terminology stems from the tradition of
having the rational map induced by the partials of $f$ called the
{\em polar map}  of the hypersurface defined by $f$.

\smallskip

Let us explain the setup of our work in a more systematic way.
Let $\Omega_{A/k}$ denote the
module of K\"ahler $k$-differentials of $A$ and let $A\simeq
k[\TT]/P$ be a presentation of $A$ over a polynomial ring $k[\TT]$.
Consider the well-known conormal exact sequence
\begin{equation}
P/P^{2}\stackrel{\delta}{\lar}\sum_{j=1}^m A\, dT_j\lar
\Omega_{A/k}\rar 0,
\end{equation}
where $\delta$ is induced by the transposed Jacobian matrix over
$k[\TT]$ of a generating set of $P$. Let $\fp\subset \sum_{j=1}^m
R\, dT_j$ denote the $R$-submodule generated by $\delta(P/P^2)$ --
the elements of which are called {\it polar syzygies\/} of $\ff$.
This module is actually a submodule of the first syzygy module $\fz$
of the transposed Jacobian module $\fd(\ff)$ when the latter is
viewed in its natural embedding in $\sum_{i=1}^nR\, dx_i$ -- the
elements of $\fz$ could be called {\it differential syzygies\/} of
$\ff$. We say that $\ff$ (or the embedding $A\subset R$) is {\it
polarizable\/} if $\fp={\cal Z}$.

\smallskip

One basic principle will tell us that, on a far more general
setting, the two modules always  have the same rank and allow for a
comparison (Lemma~\ref{ranks}).

\smallskip

When $\ff$ are monomials of degree $2$, a special case of the
presently envisaged problem had been taken up earlier in \cite{jac},
where $A$ was, up to degree normalization, the homogeneous
coordinate ring of a coordinate projection of the Segre embedding of
$\pp^r\times \pp^s$. The main result was that the $k$-subalgebra generated by
a subset of the monomials
$$\{y_i\,z_j\,|\, 0\leq i\leq r, 0\leq j\leq s\}\subset k[\,y_0,\ldots, y_r; z_0,\ldots, z_s\,]$$
is polarizable.

\smallskip

In this work we vastly enlarge the picture, obtaining a full combinatorial
characterization of polarizability.
The combinatorial gadget that plays a main role is a
graph with loops - this is allegedly a nontrivial work over the
usual simple graphs, where no loops are present.  In the more general context of admitting loops, the
given monomial generators $\ff$ of $A$ over $k$ still correspond to
(traditional) edges and loops and the corresponding graph is denoted
${\cal G}(\ff)$. Even in this generalized setting we will stick to the
terminology that has $A$ called the {\em edge-algebra} associated to
${\cal G}(\ff)$.

\smallskip

For the purpose of establishing edge-algebra polarizability, we
dwell on the fine points of the structure of both $\fp$ and $\fz$,
by describing their sets of natural minimal generators in terms of
combinatorial substructures of the corresponding graph ${\cal
G}(\ff)$. We were thus led to isolate two special configurations of
${\cal G}(\ff)$, called {\it cycle arrangements\/} and {\it
molecules}, respectively. These configurations are natural supports
of closed walks of ${\cal G}(\ff)$ and, provided these closed walks
are {\it even}, give rise to natural sets of both differential and
polar syzygies. In order to detect minimal generators among these we
further impose certain restrictions and arrive to the notion of {\it
non-split\/} and {\it indecomposable\/} even closed walks. A
consequence of these methods is a complete characterization of
polarizability in terms of the above configurations.

\smallskip

Besides throwing light into polarizability, it is to expect that
these configurations yield some new numerical invariants of the
graph that may have some curious reflection into the structure of
the corresponding algebra.

\smallskip

An almost immediate consequence is a new proof of the result that
the edge-algebra of a connected bipartite graph is polarizable --
this is precisely the main theorem in \cite[Theorem 2.3]{jac} for
the projections of the Segre embedding.

\smallskip

Using the known characterization of the integral closure of the
corresponding edge-algebra (see \cite[Theorem 1.1]{bowtie},
\cite[Corollary 2.3]{ohsugihibi}), we are able to show that
polarizability implies normality of the algebra and the converse
holds provided the graph is free of certain degenerate
configurations. Both polarizability and normality involve the
existence of the so-called {\it bow tie\/} configurations which are
special cases of the previous configurations (the terminology itself
was introduced in \cite{bowtie} and the notion was based on an
earlier construct of M. Hochster).

\smallskip

We further consider the question as to how the problem of
polarizability is affected by ``variable collapsing''  when $A$ is
generated by monomials of degree $2$. This collapsing can be thought
of as a loop-contraction operation on the edges of a graph (its
geometric interpretation in terms of Proj$(A)$ is that of projecting
down to a one dimension less ambient by cutting with a suitable
elementary hyperplane). We show that it preserves the $k$-algebra
$A$ by a $k$-isomorphism if and only the given graph is bipartite,
which can be viewed as yet another characterization of connected
bipartite graphs. Conversely, by ``resolving'' a loop issuing from
an odd cycle we improve the chances of the given generators become
polarizable.

\smallskip

From a close scrutiny of the data in a long list of computed
examples, we are naturally led to guess that there is a strong
relationship between the syzygies of the given $k$-algebra
generators $\ff$ of $A$ and polarizability. In this vein,  we first
show that the condition that the module of syzygies of $\ff$ is
generated by linear relations is equivalent to the edge graph of
${\cal G}(\ff)$ having diameter at most $2$, an easy result that
gives an algebraic tint to the notion of diameter - one would be
tempted to ask whether  the exact value of the diameter reflects a
numerical algebraic invariant, such as the dimension of the subspace
of syzygies spanned in degree $2$ (or $4$ by considering the usual
degree shift). Merging with the aforementioned combinatorial
characterization of polarizability we show that linear presentation
implies polarizability.

\smallskip

A curious consequence of the theory is that the rational map
$\pp^{n-1}\dasharrow \pp^{m-1}$ defined by a polarizable set $\ff$
of monomials of degree $2$, such that $\dim k[\ff]=n$, maps
$\pp^{n-1}$ birationally onto its image. This includes rational maps
defined by polimatroidal sets of monomials of degree $2$ of maximal
rank - a subclass of which are the so-called algebras of Veronese
type. This result recovers a couple of theorems proved in
\cite{SimisVilla} with a different approach.

\smallskip

As a final note, the reason to tackle solely monomials of degree $2$
-- and not more general toric algebras as would be the case -- is
due to an as yet not completely understood phenomenon by which such
monomial $k$-subalgebras generated in degree higher than $2$ easily
fail to be polarizable.

\section{Statement of the problem}\label{terminology}

Let $A=k[\ff]=k[f_1,\ldots, f_m]\subset R=k[x_1,\ldots, x_n]$.
Consider a presentation $A\simeq S/P$ via
$S=k[T_1,\ldots,T_m]\surjects A$ by mapping $T_j\mapsto f_j$. We
assume throughout that ${\rm char}(k)=0$.

Recall the well-known conormal sequence
\begin{equation}\label{conormal_sequence}
P/P^2\stackrel{\delta}{\lar}\sum_{j=1}^m A\, dT_j\lar \Omega_{A/k}\rar 0,
\end{equation}
where $\delta$ is induced by the transposed Jacobian matrix over $S$ of a generating
set of $P$, namely
$$\delta: F \,(\bmod {P^2})\mapsto \sum_j \frac{\partial F}
{\partial T_j} \,(\bmod P)\,dT_j.$$
The embedding $A\subset R$ induces an embedding $\sum_{j=1}^m A\, dT_j\subset
\sum_{i=1}^m R\, dT_j$.

Throughout, we set $\fp=\delta(P/P^2)\,R\subset \sum_{j=1}^m R\,
dT_j$, the $R$-submodule generated by the image of $\delta$. Then
$\fp$ is generated by the vectors $\sum_j \frac{\partial F}
{\partial T_j} \,(\ff)\,dT_j$, where $F$ runs through a set of
generators of $P$. On the other hand, by the usual rules of
composite derivatives, if $F\in P$ then $\sum_{j=1}^m \frac{\partial
F}{\partial T_j} (\ff)\, df_j=0$. This means that $\fp\subset \fz$,
where $\fz$ is the first syzygy module of the differentials $d\ff$.

\begin{Definition}\rm As a way of terminology, the elements of $\fz$
(respectively, $\fp$) are called {\it differential syzygies\/}
(respectively, {\it polar syzygies\/}). Thus, $\fz$ (respectively,
$\fp$) will be referred to as the {\it differential syzygy module\/}
(respectively, the {\it polar syzygy module\/}) of $\ff$.

The set $\ff$ (or, by a slight abuse, the embedding $A\subset R$
defined by these generators) is said to be {\it polarizable\/} if
$\fp=\fz$.
\end{Definition}

A preliminary fact in this framework is the following result, which
seems to be partially folklore (but see \cite[Proposition
1.1]{SimisDiff} for a proof and a feeling of this result
and its previous history).

\begin{Proposition}\label{dim_and_rank}
If {\rm char$(k)=0$} then $\dim k[\ff]=\rk \,\fd(\ff)$.
\end{Proposition}

It will be used in the proof of the main supporting evidence for the potential
equality $\fp=\fz$, as given by the following result of general nature.

\begin{Lemma}\label{ranks} $\rk_R(\fp)=\rk_R(\fz)\,(\,=\hht P\,)$.
\end{Lemma}
\demo Since $P/P^{(2)}$ and $\fp$ are generated by the same
generating set,  computing rank by the familiar determinantal method
yields  $\rk_A(P/P^{(2)})=\rk_R(\fp)$. But $P^{(2)}/P^2$ is a
torsion $A$-module and $P$ is generically a complete intersection on
$S$, hence $\rk_A(P/P^{(2)})=\rk_A(P/P^{2})=\hht (P)$. On the other
hand, $\rk_R({\cal Z})=m-\rk_R(\fd(\ff))=m-\dim A=\hht(P)$ using
Proposition~\ref{dim_and_rank}. Since $\fp\subset {\cal Z}$, we are
through.
\qed

\medskip

 As it turns the
theory in the case of monomials of degree $2$ is fairly under grasp;
in particular, we will
give a complete characterization of when $\ff$ is polarizable in
terms of its underlying combinatorial nature. For this, we are led
to introduce several configurations of that nature drawing largely
from the theory of graphs.

\section{Related graph substructures}

In this section we develop the graph-theoretic material needed to
translate the stated problem into combinatorics. The general
reference for algebraic graph theory in this section is
\cite{VillaBook}.

\subsection{Non-split even closed walks}

Recall that, given a set $\ff=\{f_1,\ldots,f_m\}\subset
R=k[x_1,\ldots,x_n]$ of distinct monomials of degree $2$, one
associates to it a graph ${\cal G}(\ff)$ with loops whose vertices
correspond to the variables, and where, given $i,j$, $1\leq i\leq
j\leq n$, the vertices $x_i$ and $x_j$ of ${\cal G}(\ff)$ are
connected by an edge whenever $x_ix_j\in\ff$. The ideal $(\ff)$ is
radical if and only if ${\cal G}(\ff)$ is a simple graph, i.e., has
no loops.

\medskip

The notion of even closed walk on ${\cal G}(\ff)$ is central in this
part, so let us recall its main features along with some extra
precision needed for the purpose of this paper.

An {\it even closed walk of length\/} $2r$ in ${\cal G}(\ff)$ is
given by a sequence ${\mathfrak w}=\{g_1,\ldots,g_{2r}\}$, where
$g_j\in \ff$ and $\gcd (g_{j},g_{j+1})\neq 1$ for $1\leq j\leq 2r $
(with the proviso $g_{2r+1}=g_1$). We call ${\mathfrak
w}=\{g_1,\ldots,g_{2r}\}$ the structural {\it edge sequence\/} of
the even closed walk.

Often, by abuse, we make no distinction between an even closed walk
and its structural edge sequence. Note that for a given $j$ the
corresponding edge $g_j$ may be repeated in the sequence -- we then
speak of an {\it edge repetition}. In this vein, for any edge $f$ of
the graph there is the trivial even closed walk $\{f,f\}$ -
actually, this is the only even closed walk of length $2$ in a
graph. In particular, by swinging back and forth arbitrarily often
one finds even closed walks of arbitrary length! Thus, a procedure
is needed that overlooks such useless nuisances that may creep in as
an argument gets more intricate. Such a procedure will be given soon
below.

\medskip

Note that the even closed walk ${\mathfrak w}$ may also be given by
its {\it vertex sequence}, namely:
$$g_1=x_{i_1}x_{i_2},\, g_2=x_{i_2}x_{i_3},\,\ldots,\,
g_{2r}=x_{i_{2r}}x_{i_1},$$ with $i_1,\ldots,i_{2r}\in
\{1,\ldots,n\}$. Similarly, one may have a {\it vertex repetition}.
Those even closed walks with no vertex repetition are called {\it
even cycles}. Clearly, an edge repetition implies a vertex
repetition, but not vice-versa as the following simple example
illustrates:

\begin{figure}[h]
\begin{center}
\setlength{\unitlength}{.035cm}
\begin{picture}(-60,60)(0,10)
 \put(0,0){\circle*{4}}
 \put(0,60){\circle*{4}}
 \put(-60,60){\circle*{4}}
 \put(-60,0){\circle*{4}}
 \put(-30,30){\circle*{4}}
 \put(0,60){\line(0,-1){60}}
 \put(-60,60){\line(0,-1){60}}
 \put(-30,30){\line(1,-1){30}}
 \put(-60,0){\line(1,1){30}}
 \put(0,60){\line(-1,-1){30}}
 \put(-30,30){\line(-1,1){30}}
\end{picture}
\end{center}
\caption{Path-degenerate bow tie}\label{figure_degenerate_bowtie}
\end{figure}
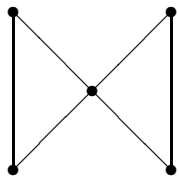

\begin{Remark}\label{uniqueness_of_edge_sequence}\rm
Let $\{g_1,\ldots,g_{2r}\}$ be the structural edge sequence of an
even closed walk ${\mathfrak w}$.
\begin{enumerate}
\item
For all $i$, $1\leq i< 2r$, the sequence
$\{g_{i+1},\ldots,g_{2r},g_1,\ldots,g_i\}$, obtained by cyclically
permuting the edges of the original edge sequence, defines the same
even closed walk ${\mathfrak w}$ as before. Therefore, by suitably
reordering the elements in an edge sequence defining ${\mathfrak
w}$, one can arbitrarily choose which of the variables
$x_{i_1},\ldots,\, x_{i_{2r}}$ comes first in a vertex sequence of
${\mathfrak w}$.
\item
There may be more ways of permuting the edges in a given edge
sequence of ${\mathfrak w}$ -- always preserving the property that
the least common multiple of two consecutive elements is not. Thus,
 in the above example of two triangles $\{f_1,f_2,f_3\}$
and $\{f_4,f_5,f_6\}$ with a common vertex belonging to $f_1$,
$f_3$, $f_4$ and $f_6$, the even closed walk with edge sequence
$\{f_1,\ldots,f_6\}$ can also be described by the sequence
$\{f_1,f_2,f_3,f_6,f_5,f_4\}$.
\end{enumerate}
\end{Remark}

Despite this lack of uniqueness of ordering of edge or vertex
sequences, we speak of them as if they were uniquely defined by the
corresponding even closed walk.

From the first of these observations follows in particular that,
given even closed walks ${\mathfrak w}_1=\{g_1,\ldots,g_{2r}\}$ and
${\mathfrak w}_2=\{g'_1,\ldots,g'_{2s}\}$ who share at least one
vertex, one can always assume that this common vertex is the first
element in their vertex sequences as observed above. We then denote
by ${\mathfrak w}_1\sqcup{\mathfrak w}_2$ the even closed walk whose
edge sequence is $\{g_1,\ldots,g_{2r},g'_1,\ldots,g'_{2s}\}$.
Conversely, one ought to consider those even closed walks that split
this way. We make this into a precise definition.

\begin{Definition}\label{defining_split}\rm
We say that an even closed walk ${\mathfrak
w}=\{g_1,\ldots,g_{2r}\}$ {\it splits\/} if it has a vertex
repetition, say the first element $x_{i_1}$ in its vertex sequence,
and if there exists $s$, $1\leq s<r$, such that ${\mathfrak
w}_1=\{g_1,\ldots,g_{2s}\}$ and ${\mathfrak
w}_2=\{g_{2s+1},\ldots,g_{2r}\}$ are the edge sequences of two
smaller even closed walks. When this occurs, we have that
${\mathfrak w}={\mathfrak w}_1\sqcup{\mathfrak w}_2$. We say that
${\mathfrak w}$ {\it splits into ${\mathfrak w}_1$ and ${\mathfrak
w}_2$}, and also that ${\mathfrak w}$ {\it splits at the vertex
$x_{i_1}$}. An even closed walk that does not split is said to be
{\it non-split}.
\end{Definition}

\begin{Remark}\rm
Even closed walk splitting has an obvious parallel in other
algebraic theories: an even closed walk containing an even closed
{\it subwalk\/} -- in the sense of a proper subset of the given edge
sequence being the edge sequence of an even closed walk -- may not
split into this and another even closed subwalk.
\end{Remark}

By definition, an even closed walk has a vertex repetition if it
splits. The converse fails as the example in
Figure~\ref{figure_degenerate_bowtie} shows. The next lemma gives
the behavior of a repeated vertex in the vertex sequence of a
non-split even closed walk. It shows in particular that the vertices
involved in a non-split even closed walk cannot occur more than
twice along its vertex sequence.

\begin{Lemma}\label{vertex_rep_in_non-split_walk}
Let ${\mathfrak w}=\{g_1,\ldots,g_{2r}\}$ be a non-split even closed
walk in ${\cal G}(\ff)$ with $g_1=x_{i_1}x_{i_2}$,
$g_2=x_{i_2}x_{i_3}$, \ldots, $g_{2r}=x_{i_{2r}}x_{i_1}$
{\rm(}$r\geq 2${\rm)}. Let $x_{i_j}$ be a repeated vertex in this
sequence, say $x_{i_j}=x_{i_l}$, with $1\leq j<l\leq 2r$. Then:
\begin{enumerate}
\item[{\rm (1)}] {\rm (Uniqueness of recurrence)}
$x_{i_k}\neq x_{i_j}$ for all $k\neq j,l$.
\item[{\rm (2)}] {\rm (Parity condition)}
$l-j\equiv 1\,(\hbox{mod}\, 2)$.
\end{enumerate}
\end{Lemma}

\demo For a suitable edge ordering, one may assume that $j=1$. Next
choose $l$ to be the smallest index such that $x_{i_l}=x_{i_1}$.
Since ${\mathfrak w}$ is non-split, $l$ has to be even, otherwise
${\mathfrak w}$ splits into $\{g_1,\ldots,g_{l-1}\}$ and
$\{g_{l},\ldots,g_{2r}\}$; hence (2) follows. Now, if
$x_{i_k}=x_{i_1}$ for some $k$, $l<k\leq 2r$, then by the same
reasoning $k$ has to be even, in which case ${\mathfrak w}$ splits
into $\{g_l,\ldots,g_{k-1}\}$ and
$\{g_k,\ldots,g_{2r},g_1,\ldots,g_{l-1}\}$; hence (1) holds as well.
\qed

\medskip
A similar result holds as regards edge repetitions in a non-split
even closed walk. Again it shows that any edge along the edge
sequence of a non-split even closed walk occurs at most twice.

\begin{Lemma}\label{edge_rep_in_non-split_walk}
Let ${\mathfrak w}=\{g_1,\ldots,g_{2r}\}$ be a non-split even closed
walk in ${\cal G}(\ff)$ with $g_1=x_{i_1}x_{i_2}$,
$g_2=x_{i_2}x_{i_3}$, \ldots, $g_{2r}=x_{i_{2r}}x_{i_1}$
{\rm(}$r\geq 2${\rm)}. If it has an edge repetition, say $g_j=g_l$
for $1\leq j<l\leq 2r$, then the following three conditions hold:
\begin{enumerate}
\item[{\rm (1)}]
{\rm (Sense-reversing recurrence)}
$x_{i_{j}}=x_{i_{l+1}}$ and $x_{i_{j+1}}=x_{i_{l}}$.
\item[{\rm (2)}]
{\rm (Uniqueness of recurrence)}
$g_k\neq g_j$ for all $k\neq j,l$.
\item[{\rm (3)}]
{\rm (Parity condition)}
$l-j\equiv 0\,(\hbox{mod}\, 2)$.
\end{enumerate}
\end{Lemma}

\demo For a suitable edge ordering, one may assume that $j=1$.

(1) Since $g_l=g_1$, one has either $x_{i_{1}}=x_{i_{l}}$ and
$x_{i_{2}}=x_{i_{l+1}}$, or $x_{i_{1}}=x_{i_{l+1}}$ and
$x_{i_{2}}=x_{i_{l}}$. If $x_{i_{1}}=x_{i_{l}}$ and
$x_{i_{2}}=x_{i_{l+1}}$, note that the edge sequence
$\{g_2,\ldots,g_{l-1},g_{2r},g_{2r-1},\ldots,g_{l+1},g_l,g_1\}$ also
defines the even closed walk ${\mathfrak w}$,  hence
 $${\mathfrak w}=\{g_2,\ldots,g_{l-1},g_{2r},g_{2r-1},\ldots,g_{l+1}\}\sqcup\{g_l,g_1\}.$$

(2) It follows from Lemma~\ref{vertex_rep_in_non-split_walk}, (1).

(3) This is clear since if $l=2s$ for some $s$, $1<s< r$, then
${\mathfrak w}$ splits into $\{g_1,\ldots, g_{2s}\}$ and
$\{g_{2s+1},\ldots,g_{2r}\}$. \qed

\subsection{Supporting configurations}

Associated to an even closed walk ${\mathfrak w}$ in ${\cal G}(\ff)$
there is a connected subgraph of ${\cal G}(\ff)$ whose edges are the
distinct elements in the edge sequence of ${\mathfrak w}$. This
subgraph of ${\cal G}(\ff)$ is called the {\it support} of
${\mathfrak w}$. Clearly, an even cycle in ${\cal G}(\ff)$ is
exactly the configuration that supports a non-split even closed walk
with no vertex repetition. As a rule, we make no distinction
between an even cycle and its naturally associated non-split even
closed walk and,  by the same abuse, we will
identify an even closed walk with its support.
We now proceed to survey a few more configurations that
support non-split even closed walks.

\subsubsection{Bow ties}

The following configuration was introduced in \cite{bowtie}.

\begin{Definition}\label{defining_bowtie}\rm
(1) A {\it bow tie\/} of ${\cal G}(\ff)$ is the (connected) subgraph
${\cal B}$ of ${\cal G}(\ff)$ consisting of two odd cycles whose
sets of edges are disjoint, connected by a unique non-empty path.
One allows for either cycle to degenerate into a loop -- in this
case, we speak of a {\it looped bow tie}.

(2) One allows the connecting path to be formed by one single edge
-- in which case we call the configuration a {\it monedge bow tie\/}
-- or to degenerate into a single vertex -- in which case we refer
to the bow tie as being {\it path-degenerate\/} (see
Figure~\ref{figure_degenerate_bowtie}). Note that, in particular, a
looped bow tie can also be a monedge (looped) bow tie, with either
or both cycles being loops; similarly, a looped bow tie can be a
path-degenerate (looped) bow tie with one of the cycles (but not
both, of course) being a loop.
\end{Definition}

Next are depicted some such configurations.

\begin{figure}[h]
\begin{center}
\setlength{\unitlength}{.035cm}
\begin{picture}(0,30)(0,0)
 \put(0,0){\circle*{4}}
 \put(-40,0){\circle*{4}}
 \put(40,0){\circle*{4}}
 \put(-60,30){\circle*{4}}
 \put(-60,-30){\circle*{4}}
 \put(-100,20){\circle*{4}}
 \put(-100,-20){\circle*{4}}
 \put(80,30){\circle*{4}}
 \put(80,-30){\circle*{4}}
 \put(-100,20){\line(0,-11){40}}
 \put(-100,20){\line(4,1){40}}
 \put(-100,-20){\line(4,-1){40}}
 \put(-40,0){\line(-2,3){20}}
 \put(-40,0){\line(-2,-3){20}}
 \put(-40,0){\line(1,0){40}}
 \put(1,0){\line(1,0){40}}
 \put(40,0){\line(4,3){40}}
 \put(40,0){\line(4,-3){40}}
 \put(80,30){\line(0,-1){60}}
\end{picture}
\end{center}
\vskip 0.7cm \caption{Typical bow tie}\label{figure_bowtie}
\end{figure}
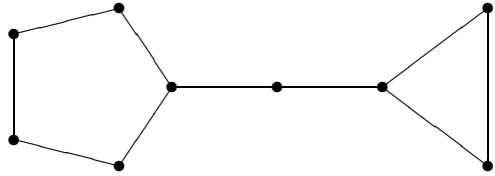

\begin{figure}[h]
\begin{center}
\setlength{\unitlength}{.035cm}
\begin{picture}(20,30)(0,0)
 \put(-130,0){\circle*{4}}
 \put(-170,30){\circle*{4}}
 \put(-170,-30){\circle*{4}}
 \put(-90,0){\circle*{4}}
 \put(-50,30){\circle*{4}}
 \put(-50,-30){\circle*{4}}
 \put(-130,0){\line(-4,3){40}}
 \put(-130,0){\line(-4,-3){40}}
 \put(-170,30){\line(0,-1){60}}
 \put(-130,0){\line(1,0){40}}
 \put(-90,0){\line(4,3){40}}
 \put(-90,0){\line(4,-3){40}}
 \put(-50,30){\line(0,-1){60}}
 \put(80,0){\circle*{4}}
 \put(40,0){\circle*{4}}
 \put(20,30){\circle*{4}}
 \put(20,-30){\circle*{4}}
 \put(-20,20){\circle*{4}}
 \put(-20,-20){\circle*{4}}
 \put(90,0){\circle{20}}
 \put(-20,20){\line(0,-1){40}}
 \put(-20,20){\line(4,1){40}}
 \put(-20,-20){\line(4,-1){40}}
 \put(40,0){\line(-2,3){20}}
 \put(40,0){\line(-2,-3){20}}
 \put(40,0){\line(1,0){40}}
 \put(140,0){\circle{20}}
 \put(150,0){\circle*{4}}
 \put(190,30){\circle*{4}}
 \put(190,-30){\circle*{4}}
 \put(150,0){\line(4,3){40}}
 \put(150,0){\line(4,-3){40}}
 \put(190,30){\line(0,-1){60}}
\end{picture}
\end{center}
\vskip 0.7cm \caption{Monedge, monedge looped and path-degenerate
looped bow ties}\label{figure_special_bowties}
\end{figure}
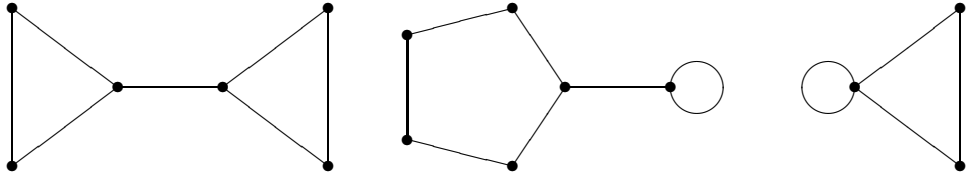

\vspace{.5cm}

Note that a bow tie is the support of a non-split even closed walk with
a vertex repetition -- indeed, an edge repetition unless it is
path-degenerate.

\medskip

\begin{Remark}\label{bowtie_and_normality}\rm These configurations were
introduced in \cite{bowtie} in order to build the integral closure
of the algebra $k[\ff]$ in case ${\cal G}(\ff)$ had no loops. The
{\it Hochster monomial\/} associated to a bow tie is the product of
the variables corresponding to the totality of the vertices of the
two cycles. We give the notion some flexibility in the sense that we
do not a priori require ${\cal B}$ to be an {\it induced\/}
subgraph, i.e., ${\cal G}(\ff)$ may have edges that do not belong to
the bow tie configuration and that connect the two structural odd
cycles, or one vertex on one odd cycle to one vertex on the path, or
one vertex on the path to another vertex on the path. The two
approaches differ in that by taking the induced subgraph definition
the Hochster monomial is a fresh generator of the integral closure
of $k[\ff]$, while our present notion allows for the Hochster
monomial to belong to $k[\ff]$ (i.e., to be a product of edges).
Otherwise, the notion is the same as in \cite{bowtie}. We will have
more to say on this theme later.
\end{Remark}

We now introduce two basic configurations in a graph which will play
a central role in this part, provided they support even closed
walks. The first one includes path-degenerate bow ties, while the
second of these configurations will be a generalized version of a
bow tie which is not path-degenerate.

\subsubsection{Cycle arrangements}

The following configuration can be thought of as an extension of the
notion of a cycle in a graph.

\begin{Definition}\rm
A {\it cycle arrangement\/} of a graph ${\cal G}(\ff)$ is a
connected subgraph of ${\cal G}(\ff)$ consisting of a set of (even
or odd) cycles, here called the {\it constituent cycles\/} of the
cycle arrangement, satisfying the following properties:
\begin{enumerate}
\item[{\rm (C$_1$)}] Any two constituent cycles have mutually disjoint edges;
\item[{\rm (C$_2$)}] Any two constituent cycles share at most one vertex;
\item[{\rm (C$_3$)}] Any vertex of the configuration belongs to at most two constituent
cycles.
\end{enumerate}
\end{Definition}

\begin{figure}[h]
\begin{center}
\setlength{\unitlength}{.035cm}
\begin{picture}(150,65)(20,10)
\put(60,70){\circle{20}} \put(30,30){\circle*{4}}
\put(30,30){\line(1,1){30}} \put(30,30){\line(1,-1){30}}
\put(60,60){\circle*{4}} \put(60,0){\circle*{4}}
\put(90,30){\circle*{4}} \put(90,30){\line(-1,1){30}}
\put(90,30){\line(-1,-1){30}} \put(90,30){\line(1,2){15}}
\put(90,30){\line(1,-2){15}} \put(105,60){\circle*{4}}
\put(105,0){\circle*{4}} \put(135,60){\circle*{4}}
\put(135,0){\circle*{4}} \put(150,30){\circle*{4}}
\put(150,30){\line(-1,2){15}} \put(105,60){\line(1,0){30}}
\put(105,0){\line(1,0){30}} \put(150,30){\line(-1,-2){15}}
\put(160,30){\circle{20}}
\end{picture}
\end{center}
\vskip -0.2cm \caption{Cycle arrangement}
\label{figure_cycle_arrangement}
\end{figure}
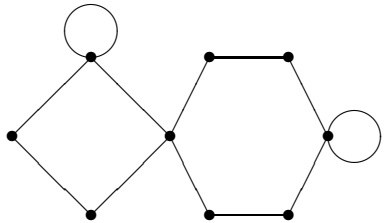

 \vskip 0.3cm \noindent
The following information ought to be kept in mind:
\begin{itemize}
\item The vertices of a cycle arrangement belonging to only one of
its constituent cycles are called {\it simple};
\item We will say that a
cycle arrangement is {\it even\/} or {\it odd\/} according to
whether the total number of edges in its configuration is even or
odd, respectively;
\item An even cycle arrangement supports an even closed walk. As often done, by abuse,
we will also refer to this even closed walk as an
even cycle arrangement. For example, the cycle arrangement in
Figure~\ref{figure_cycle_arrangement} is a non-split even closed
walk;
\item An even cycle arrangement has the property that its simple vertices
are exactly the non-repeated vertices along its vertex sequence;
\item The non-simple vertices of a cycle arrangement
belong to exactly two constituent cycles by {\rm (C$_3$)}, and hence
all vertex repetitions in a cycle arrangement satisfy the recurrence
condition Lemma~\ref{vertex_rep_in_non-split_walk}, (1).
\end{itemize}

We refer to \cite[Example~8.4.14]{VillaBook} for an example of an
even cycle arrangement which gives rise to a non-superfluous
polynomial relation of the corresponding edge algebra - we will have
more to say later about this sort of matter.

\begin{Remark}\label{splitting_cycle_arrangement}\rm
An even cycle arrangement may split. It is clear that this happens
whenever the cycle arrangement branches out into two even cycle
arrangements as shown in the following two examples:

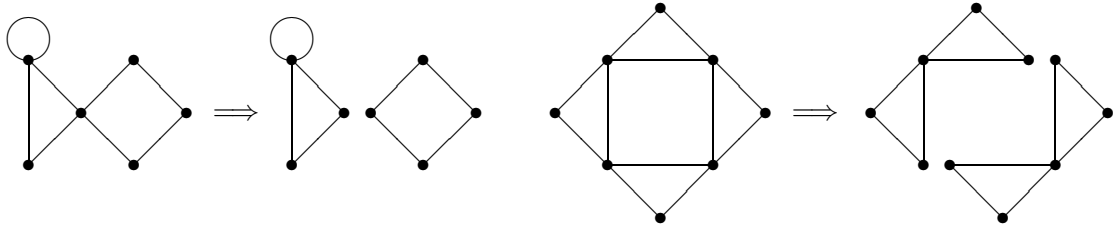
\begin{figure}[h]
\begin{center}
 \setlength{\unitlength}{.035cm}
\begin{picture}(410,70)(0,-10)
\put(0,48){\circle{16}} \put(0,40){\circle*{4}}
\put(0,40){\line(0,-1){40}} \put(0,0){\circle*{4}}
\put(0,40){\line(1,-1){20}} \put(20,20){\circle*{4}}
\put(0,0){\line(1,1){20}} \put(20,20){\line(1,-1){20}}
\put(20,20){\line(1,1){20}} \put(40,40){\circle*{4}}
\put(40,0){\circle*{4}} \put(40,40){\line(1,-1){20}}
\put(40,0){\line(1,1){20}} \put(60,20){\circle*{4}}
\put(70,17){$\mathbf \Longrightarrow$} \put(100,48){\circle{16}}
\put(100,40){\circle*{4}} \put(100,40){\line(0,-1){40}}
\put(100,0){\circle*{4}} \put(100,40){\line(1,-1){20}}
\put(120,20){\circle*{4}} \put(100,0){\line(1,1){20}}
\put(130,20){\circle*{4}} \put(130,20){\line(1,-1){20}}
\put(130,20){\line(1,1){20}} \put(150,40){\circle*{4}}
\put(150,0){\circle*{4}} \put(150,40){\line(1,-1){20}}
\put(150,0){\line(1,1){20}} \put(170,20){\circle*{4}}
\put(220,40){\circle*{4}} \put(220,0){\circle*{4}}
\put(220,40){\line(0,-1){40}} \put(220,40){\line(-1,-1){20}}
\put(220,0){\line(-1,1){20}} \put(200,20){\circle*{4}}
\put(220,40){\line(1,0){40}} \put(220,40){\line(1,1){20}}
\put(260,40){\line(-1,1){20}} \put(240,60){\circle*{4}}
\put(260,40){\circle*{4}} \put(260,0){\circle*{4}}
\put(260,40){\line(0,-1){40}} \put(260,40){\line(1,-1){20}}
\put(260,0){\line(1,1){20}} \put(280,20){\circle*{4}}
\put(220,0){\line(1,0){40}} \put(220,0){\line(1,-1){20}}
\put(260,0){\line(-1,-1){20}} \put(240,-20){\circle*{4}}
\put(290,17){$\mathbf \Longrightarrow$} \put(340,40){\circle*{4}}
\put(340,0){\circle*{4}} \put(340,40){\line(0,-1){40}}
\put(340,40){\line(-1,-1){20}} \put(340,0){\line(-1,1){20}}
\put(320,20){\circle*{4}} \put(340,40){\line(1,0){40}}
\put(340,40){\line(1,1){20}} \put(380,40){\line(-1,1){20}}
\put(360,60){\circle*{4}} \put(380,40){\circle*{4}}
\put(390,40){\circle*{4}} \put(390,0){\circle*{4}}
\put(390,40){\line(0,-1){40}} \put(390,40){\line(1,-1){20}}
\put(390,0){\line(1,1){20}} \put(410,20){\circle*{4}}
\put(350,0){\line(1,0){40}} \put(350,0){\line(1,-1){20}}
\put(390,0){\line(-1,-1){20}} \put(370,-20){\circle*{4}}
\put(350,0){\circle*{4}}
\end{picture}
\end{center}
\vskip -0.2cm \caption{Even cycle arrangements that split}
\label{figure_splitting_cycle_arrangement}
\end{figure}

 \noindent
This sort of operation will be made clear later. In
Lemma~\ref{lemma_on_nonsplit_cycle_arrangements} a characterization
will be given of when an even cycle arrangement is non-split. The
first example in Figure~\ref{figure_splitting_cycle_arrangement}
illustrates a trivial obstruction for an even cycle arrangement to
be non-split: {\it a constituent cycle that is connected to exactly
one other constituent cycle must be odd}. The second example of
Figure~\ref{figure_splitting_cycle_arrangement} puts in evidence yet
another obstruction for an even cycle arrangement to be non-split:
{\it the constituent cycles must be the only cycles of the
arrangement as a subgraph} -- thus, the inner square and the outer
octagon are not constituent cycles though they are cycles of the
containing graph. This latter necessary condition, which is however
non-obvious, will be proved in
Corollary~\ref{cycles_in_nonsplit_cycle_arrangements}.
\end{Remark}

\subsubsection{Molecules}

\medskip
We now introduce the second configuration. Recall that a {\it
path\/} of a graph is a non-closed walk without vertex repetition.
The first and last vertices of a path are called {\it extremal}.

\begin{Definition}\rm
A {\it molecule\/} of a graph ${\cal G}(\ff)$ is a connected
subgraph of ${\cal G}(\ff)$ consisting of a set of $r$ cycle
arrangements ($r\geq 2$) -- its {\it structural cycle
arrangements\/} -- and a set of $r-1$ paths -- its {\it structural
paths\/} -- satisfying the following properties:
\begin{enumerate}
\item[{\rm (M$_1$)}]
Any two structural cycle arrangements have mutually disjoint edges;
\item[{\rm (M$_2$)}]
Any two structural paths have mutually disjoint vertices (hence
mutually disjoint edges as well);
\item[{\rm (M$_3$)}]
A structural cycle arrangement and a structural path have at most
one vertex in common (in particular, have no common edges);
\item[{\rm (M$_4$)}]
Every structural cycle arrangement meets at least one structural
path and every structural path meets exactly two structural cycle
arrangements;
\item[{\rm (M$_5$)}]
Every vertex of the configuration belongs to at most two structural
cycle arrangements;
\item[{\rm (M$_6$)}]
A vertex that belongs to two structural cycle arrangements is a
simple vertex of both and a vertex that belongs to a structural
cycle arrangement and a structural path is a simple vertex of the
first and an extremal vertex of the second.
\end{enumerate}
\end{Definition}

If one draws schematically a circle for each structural cycle
arrangement and a line for each structural path, then the {\it
shadow\/} of a typical molecule is a tree (because the number of
structural paths is, by definition, one less than the number of
structural cycle arrangements), as depicted in the following
diagram:

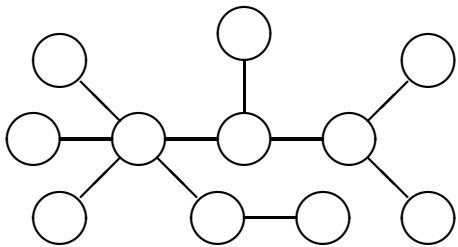
\begin{figure}[h]
\begin{center}
\setlength{\unitlength}{.035cm} \thicklines
\begin{picture}(130,85)(0,10)
\put(60,80){\circle{20}}\put(60,40){\circle{20}}
\put(60,50){\line(0,1){20}} \put(50,40){\line(-1,0){20}}
\put(20,40){\circle{20}}
\put(10,40){\line(-1,0){20}}\put(-20,40){\circle{20}}
\put(13,47){\line(-1,1){15}} \put(13,33){\line(-1,-1){15}}
\put(-10,70){\circle{20}} \put(-10,10){\circle{20}}
\put(50,10){\circle{20}} \put(27,33){\line(1,-1){15}}
\put(90,10){\circle{20}} \put(60,10){\line(1,0){20}}
\put(100,40){\circle{20}} \put(70,40){\line(1,0){20}}
\put(130,70){\circle{20}} \put(130,10){\circle{20}}
\put(107,47){\line(1,1){15}} \put(107,33){\line(1,-1){15}}
\end{picture}
\end{center}
\caption{Shadow of molecule}
\label{figure_molecule}
\end{figure}

As in the case of a cycle arrangement, the following basic
information on molecules ought to be kept in mind:
\begin{itemize}
\item
The constituent cycles of the structural cycle arrangements of a
molecule are simply called its {\it constituent cycles\/};
\item
A molecule is said to be  {\it even\/} or {\it odd\/} according as
to whether the number of edges in all its constituent cycles is even
or odd, respectively;
\item
An even molecule is the support of an even closed walk, of which all
edge repetitions correspond to the edges in the structural paths.
Again, we identify an even molecule and the even closed walk
supported on it;
\item
Each vertex in the vertex sequence of a molecule belongs to either:
one single cycle;  exactly two cycles; one single cycle and one
single path; or one single path;
\item
The vertex repetitions along the vertex
sequence of an even molecule satisfy the recurrence property
Lemma~\ref{vertex_rep_in_non-split_walk}, (1), and all its edge
repetitions satisfy the recurrence and sense-reversing properties
Lemma~\ref{edge_rep_in_non-split_walk}, (1), (2).
\end{itemize}

The previous bow tie configuration (see, e.g.,
Figure~\ref{figure_bowtie}) is a molecule with two structural cycle
arrangements consisting each of one single odd cycle (or loop), and
a single path that connects these two cycle arrangements -- the only
exception is a path-degenerate bow tie, which is a cycle arrangement
(see Definition~\ref{defining_bowtie} and the comments at the end of
the paragraph).

\medskip

\subsubsection{Skeletons of cycle arrangements and molecules}

Next one characterizes when even cycle arrangements and even
molecules are non-split. For this purpose one introduces the
following notion:

\begin{Definition}\rm
Let ${\cal B}$ be either an even cycle arrangement or an even
molecule of a graph ${\cal G}(\ff)$. The {\it skeleton\/} ${\cal
T}({\cal B})$ of ${\cal B}$ is a connected graph whose vertices fall
under two disjoint sets, the one of the {\it black\/} vertices and
the one of the {\it white\/} vertices (represented respectively by
dots and circles), defined as follows:
\begin{enumerate}
\item[{\rm (S$_1$)}]
To every constituent cycle of ${\cal B}$ there corresponds a vertex
of ${\cal T}({\cal B})$ and this vertex is black (respectively
white) if the cycle is odd (respectively, even);
\item[{\rm (S$_2$)}]
If ${\cal B}$ is a molecule then to every edge of a structural path
of ${\cal B}$ there corresponds a white vertex of ${\cal T}({\cal
B})$;
\item[{\rm (S$_3$)}]
Two vertices of ${\cal T}({\cal B})$ are connected by an edge if and
only if the corresponding sub-configurations of ${\cal B}$ --
whether constituent cycles or edges in a structural path --  meet.
\end{enumerate}
Note that the constituent cycles and the structural paths of ${\cal
B}$ uniquely determine ${\cal T}({\cal B})$.
\end{Definition}

The set of repeated vertices of ${\cal B}$ is in bijection with the
set of edges of ${\cal T}({\cal B})$. Moreover, ${\cal T}({\cal B})$
has always an even number of black vertices.

\medskip
The figures below depict the skeletons of some of the earlier
configurations. The first one is the skeleton of any monedge bow
tie. The second is the skeleton of the non-split cycle arrangement
in Figure~\ref{figure_cycle_arrangement}. The last two are the
skeletons of the two split even cycle arrangements  in
Figure~\ref{figure_splitting_cycle_arrangement}.
\begin{figure}[h]
\begin{center}
\vskip 0.5cm \setlength{\unitlength}{.035cm} \thicklines
\begin{picture}(400,20)(0,0)
\put(0,15){\circle*{6}} \put(30,15){\circle{6}}
\put(60,15){\circle*{6}} \put(0,15){\line(1,0){27}}
\put(33,15){\line(1,0){27}} \put(110,15){\circle*{6}}
\put(140,15){\circle{6}} \put(170,15){\circle{6}}
\put(200,15){\circle*{6}} \put(110,15){\line(1,0){27}}
\put(143,15){\line(1,0){24}} \put(173,15){\line(1,0){27}}
\put(250,15){\circle*{6}} \put(280,15){\circle*{6}}
\put(310,15){\circle{6}} \put(250,15){\line(1,0){57}}
\put(360,15){\circle*{6}} \put(380,35){\circle*{6}}
\put(380,-5){\circle*{6}} \put(400,15){\circle*{6}}
\put(360,15){\line(1,1){20}} \put(360,15){\line(1,-1){20}}
\put(400,15){\line(-1,1){20}} \put(400,15){\line(-1,-1){20}}
\end{picture}
\end{center}
\vskip -0.2cm \caption{Skeleton of an even molecule and of even
cycle arrangements} \label{figure_skeleton_of_cycle_arrangement}
\end{figure}
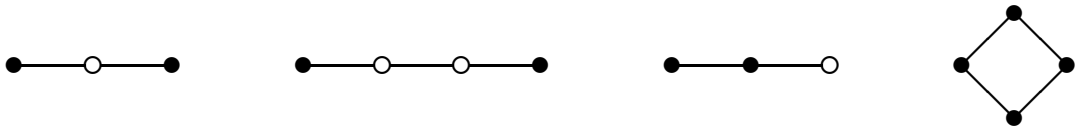

Note that if ${\cal B}$ is a molecule, its even constituent cycles
and the edges in its structural paths are represented in the same
way in ${\cal T}({\cal B})$, namely, by a white vertex. This is because an
edge in a structural path of a molecule can be considered as a
degenerate even cycle with two vertices and two edges that coincide.

\medskip
Our next result characterizes non-split even cycle arrangements and
even molecules in terms of its skeleton.

\begin{Lemma}\label{lemma_on_nonsplit_cycle_arrangements}
Let ${\cal B}$ be either an even cycle arrangement or an even
molecule of a graph ${\cal G}(\ff)$, and let ${\cal T}({\cal B})$ be
its skeleton. The following are equivalent:
\begin{enumerate}
\item[{\rm (1)}]
${\cal B}$ is non-split{\rm ;}
\item[{\rm (2)}]
No edge deletion from ${\cal T}({\cal B})$ gives rise to two
connected graphs with an even number of black
 vertices each{\rm
;}
\item[{\rm (3)}]
${\cal T}({\cal B})$ is a tree and any one edge deletion gives rise
to two trees with an odd number of black vertices each.
\end{enumerate}
\end{Lemma}

\demo The contrapositive of the implication (1) $\Rightarrow$ (2)
is straightforward by recalling that an
edge of ${\cal T}({\cal B})$ corresponds to a vertex repetition in
${\cal B}$, and that an even closed walk that splits will do so at
one of its vertex repetitions. Actually, the negation of (2) is a
reformulation of the phenomenon described in
Remark~\ref{splitting_cycle_arrangement}.

(3) $\Rightarrow$ (1): Assume that ${\cal T}({\cal B})$ is a tree
and that ${\cal B}$ splits. As already observed, this will happen at
one of its vertex repetitions and hence, removing the corresponding
edge of ${\cal T}({\cal B})$, one obtains two trees with an even
number of black vertices each.

(2) $\Rightarrow$ (3): We will be done if we show that (2) implies
that ${\cal T}({\cal B})$ is a tree because the second part of (3)
then trivially holds. Let us assume that ${\cal T}({\cal B})$ has at
least one cycle and show that (2) fails. This will be proved by
induction on the number of cycles of ${\cal T}({\cal B})$. If it has
one single cycle, say ${\cal T}_1$, then removing any two edges of
${\cal T}_1$, one gets two trees, and one only has to prove that one
can always choose two edges of ${\cal T}_1$ such that both trees
have an even number of black vertices. If there exists one vertex
${\bf v}$ in ${\cal T}_1$ such that, removing the two edges of
${\cal T}_1$ going through ${\bf v}$, one gets two trees with an
even number of black vertices, we are done. Otherwise, each vertex
in ${\cal T}_1$ satisfies that, removing the two edges of ${\cal
T}_1$ going through it, one gets two trees with an odd number of
black vertices. Now consider any two consecutive vertices of ${\cal
T}_1$ and remove from ${\cal T}_1$ the edge that goes through each
of them and which is distinct from the edge that connects them. One
gets the two expected trees. Finally, if ${\cal T}({\cal B})$ has
more than one cycle, note that removing one edge in one of its
cycles, one gets a connected graph with one cycle less, and by
induction we are done. \qed

\begin{Corollary}\label{cycles_in_nonsplit_cycle_arrangements}
A non-split even cycle arrangement or a non-split even molecule of a
graph ${\cal G}(\ff)$ includes no other cycle of ${\cal G}(\ff)$
other than its constituent cycles.
\end{Corollary}

\demo By the previous lemma, the skeleton of a non-split even cycle
arrangement or a non-split even molecule ${\cal B}$ is a tree, hence
there cannot be any additional cycles of ${\cal G}(\ff)$ in ${\cal
B}$ other than its constituent cycles. \qed

\medskip

A non-split even cycle arrangement has vertex repetitions (unless it
is a single cycle) and no edge repetition. A non-split even molecule
has always edge repetitions. The following result states that these
are all possible non-split even closed walks in ${\cal G}(\ff)$.

\begin{Proposition}\label{all_non-split_even_walks}
A non-split even closed walk in a graph ${\cal G}(\ff)$ is either an
even cycle arrangement or an even molecule.
\end{Proposition}

This result is a direct consequence of
Lemmas~\ref{vertex_rep_in_non-split_walk} and
\ref{edge_rep_in_non-split_walk} and the following two lemmas:

\begin{Lemma}\label{lemma_on_cycle_arrangements}
Let ${\mathfrak w}=\{g_1,\ldots,g_{t}\}$ be a closed walk {\rm(}even
or odd{\rm)} in ${\cal G}(\ff)$ with $g_1=x_{i_1}x_{i_2}$,
$g_2=x_{i_2}x_{i_3}$, \ldots, $g_{t}=x_{i_{t}}x_{i_1}$. Assume that
${\mathfrak w}$ has no edge repetition, and that any vertex
repetition $x_{i_{j}}=x_{i_{l}}$ for $1\leq j<l\leq t$ satisfies the
recurrence condition of {\rm
Lemma~\ref{vertex_rep_in_non-split_walk}}. Then, ${\mathfrak w}$ is
a cycle arrangement.
\end{Lemma}

\demo The proof is by induction on the number $s\geq 0$ of vertex
repetitions in ${\mathfrak w}$. If $s=0$, then ${\mathfrak w}$ is a
cycle. If $s\geq 1$, one can assume without loss of generality that
$x_{i_1}$ is a vertex repetition, i.e., $x_{i_1}=x_{i_l}$ for some
$l$, $1<l\leq t$, and that $x_{i_1},\ldots,x_{i_{l-1}}$ are all
distinct (there is always a vertex repetition with this property).
Then, ${\mathfrak w}':=\{g_1,\ldots,g_{l-1}\}$ is a cycle in ${\cal
G}(\ff)$ whose vertex sequence contains $x_{i_1}$, and ${\mathfrak
w}'':=\{g_l,\ldots,g_{t}\}$ is a closed walk in ${\cal G}(\ff)$
whose vertex sequence contains $x_{i_1}$ with $s-1$ vertex
repetitions ($x_{i_1}$ is not a vertex repetition in ${\mathfrak
w}''$) that satisfies the recurrence condition of
Lemma~\ref{vertex_rep_in_non-split_walk}. Applying the recursive
hypothesis we are done. \qed

\begin{Lemma}\label{lemma_on_molecules}
Let ${\mathfrak w}=\{g_1,\ldots,g_{t}\}$ be a closed walk {\rm(}even
or odd{\rm)} in ${\cal G}(\ff)$ with $g_1=x_{i_1}x_{i_2}$,
$g_2=x_{i_2}x_{i_3}$, \ldots, $g_{t}=x_{i_{t}}x_{i_1}$ satisfying
that $g_j\neq g_{j+1}$ for all $j=1,\ldots,t$ {\rm(}with the proviso
$g_{t+1}=g_1${\rm)}. Assume that any vertex repetition
$x_{i_{j}}=x_{i_{l}}$ for $1\leq j<l\leq t$ satisfies the recurrence
condition of $\,${\rm Lemma~\ref{vertex_rep_in_non-split_walk}}, and
that any edge repetition $g_j=g_l$ for $1\leq j<l\leq t$ satisfies
the sense-reversing property in {\rm
Lemma~\ref{edge_rep_in_non-split_walk}}. Then, ${\mathfrak w}$ is a
molecule.
\end{Lemma}

\demo The proof is by induction on the number $s\geq 0$ of edge
repetitions in ${\mathfrak w}$. The case $s=0$ is
Lemma~\ref{lemma_on_cycle_arrangements}. If $s\geq 1$, one can
assume without loss of generality that $g_1$ is an edge repetition,
i.e., $g_\ell=g_1$ for some $\ell$, $2<\ell\leq t$, and that
$g_1,\ldots,g_{\ell-1}$ are all distinct (there is always at least
one edge repetition with this property). By the sense-reversing
property (1) in Lemma~\ref{edge_rep_in_non-split_walk},
$x_{i_\ell}=x_{i_2}$ and $x_{i_{\ell+1}}=x_{i_1}$, and hence
$g_{\ell-1}=x_{i_{\ell-1}}x_{i_2}$ and
$g_{\ell+1}=x_{i_1}x_{i_{\ell+2}}$. Thus, ${\mathfrak
w}':=\{g_2,\ldots,g_{\ell-1}\}$ is a closed walk in ${\cal G}(\ff)$
and since it has no edge repetition, it is a cycle arrangement by
Lemma~\ref{lemma_on_cycle_arrangements}. Set ${\mathfrak
w}'':=\{g_{\ell},\ldots,g_{t},g_1\}$. It is a closed walk in ${\cal
G}(\ff)$ whose vertex and edge repetitions satisfy the same
properties as the ones in ${\mathfrak w}$ because ${\mathfrak
w}={\mathfrak w}'\sqcup{\mathfrak w}''$. Since $g_{\ell}=g_1$, by
the sense-reversing property (1) in
Lemma~\ref{edge_rep_in_non-split_walk}, one has that
$\{g_{\ell+1},\ldots,g_{t}\}$ is also a closed walk in ${\cal
G}(\ff)$. By the same argument, if $g_{\ell+1}=g_t$, then
$\{g_{\ell+2},\ldots,g_{t-1}\}$ is a closed walk in ${\cal G}(\ff)$
and, iterating, we get that for some $k\geq 0$, ${\mathfrak
w}''':=\{g_{\ell+1+k},\ldots,g_{t-k}\}$ is a closed walk in ${\cal
G}(\ff)$ with $g_{\ell+1+k}\neq g_{t-k}$. This closed walk satisfies
the same conditions as ${\mathfrak w}$ and it has at most $s-1$ edge
repetitions (more precisely it has $s-1-k$ edge repetitions). If
$s-1-k\neq 0$, applying the recursive hypothesis to ${\mathfrak
w}'''$, one gets that it is a molecule. Otherwise, it is a cycle
arrangement by Lemma~\ref{lemma_on_cycle_arrangements}. We conclude
observing that the original configuration supporting the closed walk
${\mathfrak w}$ is exactly the one obtained by connecting the cycle
arrangement ${\mathfrak w}'$ to the molecule (or cycle arrangement
when $s-1-k=0$) ${\mathfrak w}'''$ by the path supported on
$\{g_\ell,\ldots,g_{\ell+k}\}$. By the recurrence property (1) in
Lemma~\ref{vertex_rep_in_non-split_walk}, the vertex
$x_{i_{\ell}}(=x_{i_2})$, respectively $x_{i_{\ell+k+1}}$, is not a
repetition in the vertex sequence of ${\mathfrak w}'$, respectively
${\mathfrak w}'''$, and hence ${\mathfrak w}$ is a molecule. \qed

\medskip

Thus, the non-split even closed walks in a graph ${\cal G}(\ff)$ are
exactly its non-split even cycle arrangements and its non-split even
molecules which are characterized in
Lemma~\ref{lemma_on_nonsplit_cycle_arrangements}.

\subsection{Indecomposable even closed walks}

We now introduce a subtler class of non-split even closed walks that
will tie up polarizability of $\ff$ to combinatorial properties of
the graph ${\cal G}(\ff)$.

\begin{Definition}\label{defining_indecomposable}\rm
A non-split even closed walk ${\mathfrak w}$ in a graph ${\cal
G}(\ff)$ is {\it decomposable\/} if there exist
$h_{1},\ldots,h_{t}\in\ff$ satisfying the following conditions:
\begin{enumerate}
\item[{\rm (D$_1$)}] $h_{1},\ldots,h_{t}$ are square free;
\item[{\rm (D$_2$)}] Any variable involved in the monomials $h_{1},\ldots,h_{t}$
corresponds to a vertex along the vertex sequence of ${\mathfrak w}$;
\item[{\rm (D$_3$)}] By adding twice every $h_j$ to the edge sequence of ${\mathfrak w}$
then, up to conveniently reordering the resulting sequence, one gets
an even closed walk that splits into two smaller even closed walks
${\mathfrak w}_1$ and ${\mathfrak w}_2$ that do not contain
${\mathfrak w}$ and whose edge sequences both contain
$h_{1},\ldots,h_{t}$.
\end{enumerate}
\end{Definition}

Next are a few simple examples of decomposable even closed walks to
bear in mind:

\newpage

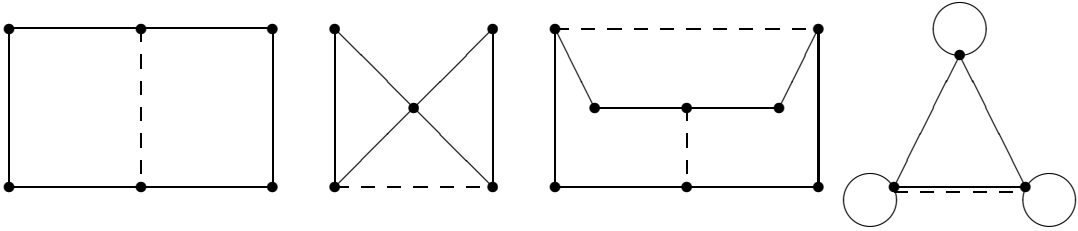
\begin{figure}[h]
\begin{center}
\setlength{\unitlength}{.035cm}
\begin{picture}(190,60)(0,10)
 \put(-100,0){\circle*{4}}
 \put(0,0){\circle*{4}}
 \put(-100,60){\circle*{4}}
 \put(0,60){\circle*{4}}
 \put(-50,0){\circle*{4}}
 \put(-50,60){\circle*{4}}
 \put(-50,60){\line(0,-1){5}}
 \put(-50,50){\line(0,-1){5}}
 \put(-50,40){\line(0,-1){5}}
 \put(-50,30){\line(0,-1){5}}
 \put(-50,20){\line(0,-1){5}}
 \put(-50,10){\line(0,-1){5}}
 \put(-100,0){\framebox(100,60){}}
\end{picture}
\begin{picture}(0,60)(0,10)
 \put(-110,0){\circle*{4}}
 \put(-110,60){\circle*{4}}
 \put(-170,60){\circle*{4}}
 \put(-170,0){\circle*{4}}
 \put(-140,30){\circle*{4}}
 \put(-110,60){\line(0,-1){60}}
 \put(-170,60){\line(0,-1){60}}
 \put(-140,30){\line(1,-1){30}}
 \put(-170,0){\line(1,1){30}}
 \put(-110,60){\line(-1,-1){30}}
 \put(-140,30){\line(-1,1){30}}
 \put(-170,0){\line(1,0){5}}
 \put(-160,0){\line(1,0){5}}
 \put(-150,0){\line(1,0){5}}
 \put(-140,0){\line(1,0){5}}
 \put(-130,0){\line(1,0){5}}
 \put(-120,0){\line(1,0){5}}
\end{picture}
\begin{picture}(0,60)(0,10)
 \put(-90,0){\circle*{4}}
 \put(10,0){\circle*{4}}
 \put(-90,60){\circle*{4}}
 \put(10,60){\circle*{4}}
 \put(-40,0){\circle*{4}}
 \put(-40,30){\circle*{4}}
 \put(-75,30){\circle*{4}}
 \put(-5,30){\circle*{4}}
 \put(-40,30){\line(0,-1){5}}
 \put(-40,20){\line(0,-1){5}}
 \put(-40,10){\line(0,-1){5}}
 \put(-90,60){\line(1,-2){15}}
 \put(10,60){\line(-1,-2){15}}
 \put(-75,30){\line(1,0){70}}
 \put(-90,0){\line(0,1){60}}
 \put(10,0){\line(0,1){60}}
 \put(-90,0){\line(1,0){100}}
 \put(-90,60){\line(1,0){5}}
 \put(-80,60){\line(1,0){5}}
 \put(-70,60){\line(1,0){5}}
 \put(-60,60){\line(1,0){5}}
 \put(-50,60){\line(1,0){5}}
 \put(-40,60){\line(1,0){5}}
 \put(-30,60){\line(1,0){5}}
 \put(-20,60){\line(1,0){5}}
 \put(-10,60){\line(1,0){5}}
 \put(0,60){\line(1,0){5}}
\end{picture}
\begin{picture}(0,60)(0,10)
 \put(35,-2){\line(1,0){5}}
 \put(45,-2){\line(1,0){5}}
 \put(55,-2){\line(1,0){5}}
 \put(65,-2){\line(1,0){5}}
 \put(75,-2){\line(1,0){5}}
 \put(60,50){\circle*{4}}
 \put(60,60){\circle{20}}
 \put(35,0){\circle*{4}}
 \put(26,-5){\circle{20}}
 \put(85,0){\circle*{4}}
 \put(94,-5){\circle{20}}
 \put(35,0){\line(1,2){25}}
 \put(85,0){\line(-1,2){25}}
 \put(35,0){\line(1,0){50}}
\end{picture}
\end{center}
\caption{Decomposable even closed walks}\label{figure_decomposable}
\end{figure}

We say that the set $h_{1},\ldots,h_{t}$ is a {\it decomposing
set\/} of ${\mathfrak w}$. An {\it indecomposable\/} even closed
walk is a non-split even closed walk which is not decomposable. In
the examples illustrated in Figure~\ref{figure_decomposable} the dotted
edges are the decomposing edges in each case. Note that every $h_j$
may belong to the very edge sequence of ${\mathfrak w}$ as the
fourth example in Figure~\ref{figure_decomposable} shows.

\begin{figure}[h]
\begin{center}
\setlength{\unitlength}{.035cm}
\begin{picture}(120,10)(0,20)
 \put(0,20){\circle{20}}
 \put(60,20){\circle{20}}
 \put(120,20){\circle{20}}
 \put(0,10){\line(1,0){120}}
 \put(0,10){\circle*{4}}
 \put(60,10){\circle*{4}}
 \put(120,10){\circle*{4}}
\end{picture}
\end{center}
\caption{Indecomposable looped bow tie}\label{figure_indecomposable}
\end{figure}
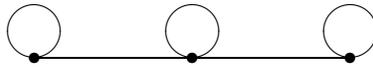

\medskip
The following example illustrates the role of condition (D$_1$) in
the definition of decomposability: in the graph in
Figure~\ref{figure_indecomposable}, if one considers the looped bow
tie involving the first and third loops, it is indecomposable since
one cannot use the second loop to decompose it because the monomial
corresponding to a loop is not square free.

\medskip

Among the non-split even closed walks in ${\cal G}(\ff)$, many are
decomposable as the following result shows:

\begin{Lemma}\label{lemma_on_number_of_cycles_in_indecomposable}
Let ${\mathfrak w}$ be a non-split even closed walk in a graph
${\cal G}(\ff)$. Assume that ${\mathfrak w}$ contains a cycle of
which at least two vertices are vertex repetitions of ${\mathfrak
w}$. Then, ${\mathfrak w}$ is decomposable.
\end{Lemma}
\demo In order to prove the result, we will show that if there
exists such a cycle ${\cal C}$, one can use the elements in $\ff$
corresponding to some of its edges as decomposing set. Note that
(D$_1$) and (D$_2$) will always be satisfied if $h_1,\ldots,h_t$
correspond to edges of ${\cal C}$, so we have to select them such
that (D$_3$) holds.

By Proposition~\ref{all_non-split_even_walks}, ${\mathfrak w}$ is
either an even cycle arrangement or an even molecule, and ${\cal C}$
is one of its constituent cycles by
Corollary~\ref{cycles_in_nonsplit_cycle_arrangements}. Our
assumption is that there are two variables, say $x_{i_1}$ and
$x_{i_2}$, corresponding to vertices along the vertex sequence of
${\cal C}$, that are vertex repetitions of ${\mathfrak w}$. By
Lemma~\ref{lemma_on_nonsplit_cycle_arrangements}~(3), removing from
the skeleton ${\cal T}({\mathfrak w})$ of ${\mathfrak w}$ the edge
corresponding to the vertex repetition $x_{i_1}$, one gets two trees
with an odd number of black vertices. One of them contains the
vertex of ${\cal T}({\mathfrak w})$ associated to the constituent
cycle ${\cal C}$ of ${\mathfrak w}$, and one does not. Denote by
${\cal G}_1$ the subgraph of ${\mathfrak w}$ corresponding to the
later. It is the support of an odd closed walk whose vertex sequence
contains $x_{i_1}$ as a non repeated vertex. We define similarly
${\cal G}_2$ by substituting $x_{i_2}$ for $x_{i_1}$. Choose any of
the two paths in ${\cal C}$ connecting $x_{i_1}$ and $x_{i_2}$, and
consider the even molecule ${\mathfrak w}_1$ obtained connecting
${\cal G}_1$ to ${\cal G}_2$ by this path. On the other hand,
consider the even closed walk ${\mathfrak w}_2$ supported by the
subgraph of ${\mathfrak w}$ obtained by removing ${\cal G}_1$ and
${\cal G}_2$. One can now easily check that (D$_3$) holds for the
decomposing set $h_{1},\ldots,h_{t}$ corresponding to the edges of
the cycle ${\cal C}$ connecting $x_{i_1}$ and $x_{i_2}$ that we have
chosen before. \qed

\medskip

The next picture illustrates the idea of the proof with an example:

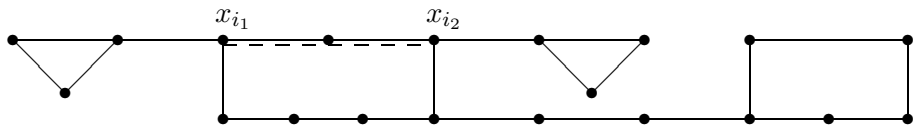
\begin{figure}[h]
\begin{center}
\vskip 0.8cm \setlength{\unitlength}{.035cm}
\begin{picture}(300,20)(0,-5)
 \put(-20,20){\circle*{4}}\put(-20,20){\line(1,0){40}}
 \put(-20,20){\line(1,-1){20}}
 \put(0,0){\circle*{4}}
 \put(20,20){\circle*{4}} \put(0,0){\line(1,1){20}}
 \put(20,20){\line(1,0){40}}
 \put(60,20){\circle*{4}} \put(100,20){\circle*{4}}
 \put(140,20){\circle*{4}}
 \put(60,20){\line(1,0){80}}
 \put(60,20){\line(0,-1){30}} \put(60,-10){\circle*{4}}
 \put(140,20){\line(0,-1){30}} \put(140,-10){\circle*{4}}
 \put(60,-10){\line(1,0){80}}
 \put(87,-10){\circle*{4}} \put(113,-10){\circle*{4}}
 \put(57,27){$x_{i_1}$} \put(137,27){$x_{i_2}$}
 \put(140,20){\line(1,0){80}}
 \put(180,20){\circle*{4}}\put(200,0){\circle*{4}}
 \put(220,20){\circle*{4}}
 \put(180,20){\line(1,-1){20}} \put(200,0){\line(1,1){20}}
 \put(140,-10){\line(1,0){120}}
 \put(180,-10){\circle*{4}} \put(220,-10){\circle*{4}}
 \put(260,20){\circle*{4}}\put(260,-10){\circle*{4}}
 \put(260,20){\line(0,-1){30}} \put(260,20){\line(1,0){60}}
 \put(320,20){\circle*{4}} \put(320,20){\line(0,-1){30}}
 \put(320,-10){\circle*{4}}
 \put(260,-10){\line(1,0){60}}
 \put(290,-10){\circle*{4}}
\put(60,18){\line(1,0){5}} \put(70,18){\line(1,0){5}}
\put(80,18){\line(1,0){5}} \put(90,18){\line(1,0){5}}
\put(100,18){\line(1,0){5}} \put(110,18){\line(1,0){5}}
\put(120,18){\line(1,0){5}} \put(130,18){\line(1,0){5}}
\end{picture}
\end{center}

\vskip -0.2cm \caption{A non-split even molecule with more than 3
constituent cycles is decomposable}
\label{figure_too_many_cycles_in_molecule}
\end{figure}

\medskip

As a consequence of
Lemma~\ref{lemma_on_number_of_cycles_in_indecomposable}, one gets
the following result which establishes a characterization of the
subclass of indecomposable even closed walks in parallel to
Proposition~\ref{all_non-split_even_walks}:

\begin{Proposition}\label{all_indecomposable_even_walks}
An indecomposable even closed walk in a graph ${\cal G}(\ff)$ is
either an even cycle or a bow tie.
\end{Proposition}

\demo An indecomposable even walk is non-split and hence, it is
either a non-split even cycle arrangement or a non-split even
molecule by Proposition~\ref{all_non-split_even_walks}. Moreover, it
cannot have more than two constituent cycles, otherwise at least one
of them would satisfy the hypothesis in
Lemma~\ref{lemma_on_number_of_cycles_in_indecomposable}. If it has
one, it is a cycle. Otherwise, it is a bow tie (that can be
path-degenerate or not). \qed

\medskip

One can now tell exactly all the indecomposable even closed walks.
We will say that a subgraph of ${\cal G}(\ff)$ is {\it induced\/} if
it is obtained by deleting a set of vertices and all the edges that
go through them and/or by deleting a set of loops (leaving of course
the base vertex of the loop). One realizes that this is the usual
concept for simple graphs, taking care in addition of loops as well.

\begin{Proposition}\label{indecomposable_are_cycles_and_induced_bowties}
The indecomposable even closed walks of a graph ${\cal G}(\ff)$ are
its indecomposable even cycles and its induced bow ties.
\end{Proposition}

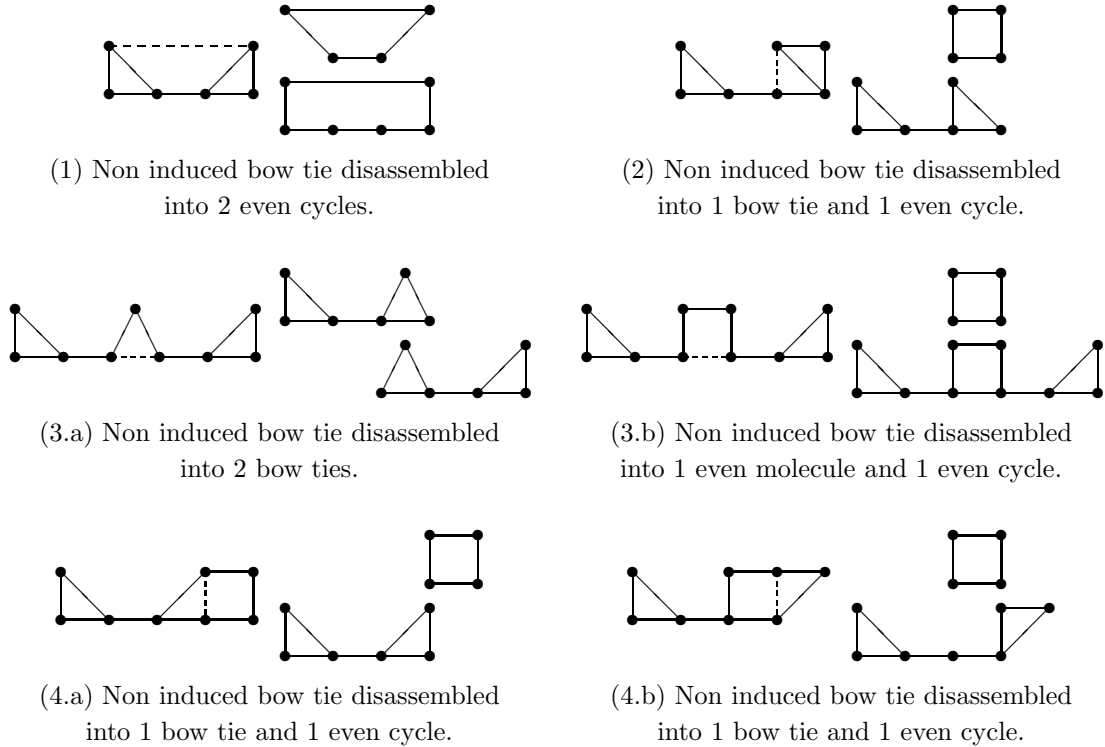
\begin{figure}[h]
\vskip 0.5cm \setlength{\unitlength}{.032cm}
\begin{center}
\begin{tabular}{cc}
\begin{tabular}{cc}
\begin{picture}(60,40)(0,-20)
 \put(0,0){\circle*{4}}
 \put(0,20){\circle*{4}}
 \put(20,0){\circle*{4}}
 \put(40,0){\circle*{4}}
 \put(60,0){\circle*{4}}
 \put(60,20){\circle*{4}}
 \put(0,0){\line(1,0){60}}
 \put(0,0){\line(0,1){20}}
 \put(0,20){\line(1,-1){20}}
 \put(40,0){\line(1,1){20}}
 \put(60,20){\line(0,-1){20}}
 \put(0,20){\line(1,0){3}} \put(6,20){\line(1,0){3}}
 \put(12,20){\line(1,0){3}} \put(18,20){\line(1,0){3}}
 \put(24,20){\line(1,0){3}} \put(30,20){\line(1,0){3}}
 \put(36,20){\line(1,0){3}} \put(42,20){\line(1,0){3}}
 \put(48,20){\line(1,0){3}} \put(54,20){\line(1,0){3}}
\end{picture}
&
\begin{picture}(60,50)(0,-20)
 \put(0,35){\circle*{4}}
 \put(20,15){\circle*{4}}
 \put(40,15){\circle*{4}}
 \put(60,35){\circle*{4}}
 \put(20,15){\line(1,0){20}}
 \put(0,35){\line(1,-1){20}}
 \put(40,15){\line(1,1){20}}
 \put(0,35){\line(1,0){60}}
 \put(0,-15){\circle*{4}}
 \put(0,5){\circle*{4}}
 \put(20,-15){\circle*{4}}
 \put(40,-15){\circle*{4}}
 \put(60,-15){\circle*{4}}
 \put(60,5){\circle*{4}}
 \put(0,-15){\line(1,0){60}}
 \put(0,-15){\line(0,1){20}}
 \put(60,5){\line(0,-1){20}}
 \put(0,5){\line(1,0){60}}
\end{picture}
\end{tabular}
&
\begin{tabular}{cc}
\begin{picture}(60,40)(0,-20)
 \put(0,0){\circle*{4}}
 \put(0,20){\circle*{4}}
 \put(20,0){\circle*{4}}
 \put(40,0){\circle*{4}}
 \put(60,0){\circle*{4}}
 \put(60,20){\circle*{4}}
 \put(40,20){\circle*{4}}
 \put(0,0){\line(1,0){60}}
 \put(0,0){\line(0,1){20}}
 \put(0,20){\line(1,-1){20}}
 \put(40,20){\line(1,-1){20}}
 \put(40,20){\line(1,0){20}}
 \put(60,20){\line(0,-1){20}}
 \put(40,0){\line(0,1){2}}\put(40,4){\line(0,1){2}}
 \put(40,8){\line(0,1){2}}\put(40,12){\line(0,1){2}}
 \put(40,16){\line(0,1){2}}
\end{picture}
&
\begin{picture}(60,50)(0,-20)
 \put(40,15){\circle*{4}}
 \put(60,15){\circle*{4}}
 \put(60,35){\circle*{4}}
 \put(40,35){\circle*{4}}
 \put(40,15){\line(1,0){20}}
 \put(40,35){\line(1,0){20}}
 \put(60,35){\line(0,-1){20}}
 \put(40,15){\line(0,1){20}}
 \put(0,-15){\circle*{4}}
 \put(0,5){\circle*{4}}
 \put(20,-15){\circle*{4}}
 \put(40,-15){\circle*{4}}
 \put(60,-15){\circle*{4}}
 \put(40,5){\circle*{4}}
 \put(0,-15){\line(1,0){60}}
 \put(0,-15){\line(0,1){20}}
 \put(0,5){\line(1,-1){20}}
 \put(40,5){\line(1,-1){20}}
 \put(40,-15){\line(0,1){20}}
\end{picture}
\end{tabular}
\\
 {\small (1)
 Non induced bow tie disassembled}
&
 {\small (2)
 Non induced bow tie disassembled}
\\
 {\small into 2 even cycles.}
&
 {\small into 1 bow tie and 1 even cycle.}
\\ \\
\begin{tabular}{cc}
\begin{picture}(99,40)(0,-20)
 \put(0,0){\circle*{4}}
 \put(0,20){\circle*{4}}
 \put(20,0){\circle*{4}}
 \put(40,0){\circle*{4}}
 \put(60,0){\circle*{4}}
 \put(80,0){\circle*{4}}
 \put(100,0){\circle*{4}}
 \put(100,20){\circle*{4}}
 \put(50,20){\circle*{4}}
 \put(0,0){\line(1,0){40}}
 \put(60,0){\line(1,0){40}}
 \put(0,0){\line(0,1){20}}
 \put(0,20){\line(1,-1){20}}
 \put(80,0){\line(1,1){20}}
 \put(100,20){\line(0,-1){20}}
 \put(40,0){\line(1,2){10}}
 \put(60,0){\line(-1,2){10}}
 \put(40,0){\line(1,0){2}} \put(44,0){\line(1,0){2}}
 \put(48,0){\line(1,0){2}} \put(52,0){\line(1,0){2}}
 \put(56,0){\line(1,0){2}}
\end{picture}
&
\begin{picture}(99,60)(0,-20)
 \put(0,15){\circle*{4}}
 \put(0,35){\circle*{4}}
 \put(20,15){\circle*{4}}
 \put(40,15){\circle*{4}}
 \put(60,15){\circle*{4}}
 \put(50,35){\circle*{4}}
 \put(0,15){\line(1,0){40}}
 \put(0,15){\line(0,1){20}}
 \put(0,35){\line(1,-1){20}}
 \put(40,15){\line(1,2){10}}
 \put(60,15){\line(-1,2){10}}
 \put(40,15){\line(1,0){20}}
 \put(40,-15){\circle*{4}}
 \put(60,-15){\circle*{4}}
 \put(80,-15){\circle*{4}}
 \put(100,-15){\circle*{4}}
 \put(100,5){\circle*{4}}
 \put(50,5){\circle*{4}}
 \put(60,-15){\line(1,0){40}}
 \put(80,-15){\line(1,1){20}}
 \put(100,5){\line(0,-1){20}}
 \put(40,-15){\line(1,2){10}}
 \put(60,-15){\line(-1,2){10}}
 \put(40,-15){\line(1,0){20}}
\end{picture}
\end{tabular}
&
\begin{tabular}{cc}
\begin{picture}(99,40)(0,-20)
 \put(0,0){\circle*{4}}
 \put(0,20){\circle*{4}}
 \put(20,0){\circle*{4}}
 \put(40,0){\circle*{4}}
 \put(60,0){\circle*{4}}
 \put(80,0){\circle*{4}}
 \put(100,0){\circle*{4}}
 \put(100,20){\circle*{4}}
 \put(40,20){\circle*{4}}
 \put(60,20){\circle*{4}}
 \put(0,0){\line(1,0){40}}
 \put(60,0){\line(1,0){40}}
 \put(0,0){\line(0,1){20}}
 \put(0,20){\line(1,-1){20}}
 \put(80,0){\line(1,1){20}}
 \put(100,20){\line(0,-1){20}}
 \put(40,0){\line(0,1){20}}
 \put(60,0){\line(0,1){20}}
 \put(40,20){\line(1,0){20}}
 \put(40,0){\line(1,0){2}} \put(44,0){\line(1,0){2}}
 \put(48,0){\line(1,0){2}} \put(52,0){\line(1,0){2}}
 \put(56,0){\line(1,0){2}}
\end{picture}
&
\begin{picture}(99,60)(0,-20)
 \put(40,15){\circle*{4}}
 \put(60,15){\circle*{4}}
 \put(40,35){\circle*{4}}
 \put(60,35){\circle*{4}}
 \put(40,15){\line(1,0){20}}
 \put(40,15){\line(0,1){20}}
 \put(60,15){\line(0,1){20}}
 \put(40,35){\line(1,0){20}}
 \put(0,-15){\circle*{4}}
 \put(0,5){\circle*{4}}
 \put(20,-15){\circle*{4}}
 \put(40,-15){\circle*{4}}
 \put(60,-15){\circle*{4}}
 \put(80,-15){\circle*{4}}
 \put(100,-15){\circle*{4}}
 \put(100,5){\circle*{4}}
 \put(40,5){\circle*{4}}
 \put(60,5){\circle*{4}}
 \put(0,-15){\line(1,0){40}}
 \put(60,-15){\line(1,0){40}}
 \put(0,-15){\line(0,1){20}}
 \put(0,5){\line(1,-1){20}}
 \put(80,-15){\line(1,1){20}}
 \put(100,5){\line(0,-1){20}}
 \put(40,-15){\line(0,1){20}}
 \put(60,-15){\line(0,1){20}}
 \put(40,5){\line(1,0){20}}
 \put(40,-15){\line(1,0){20}}
\end{picture}
\end{tabular}
\\
 {\small (3.a)
 Non induced bow tie disassembled}
&
 {\small (3.b)
 Non induced bow tie disassembled}
\\
 {\small into 2 bow ties.}
&
 {\small into 1 even molecule and 1 even cycle.}
\\ \\
\begin{tabular}{cc}
\begin{picture}(80,40)(0,-20)
 \put(0,0){\circle*{4}}
 \put(0,20){\circle*{4}}
 \put(20,0){\circle*{4}}
 \put(40,0){\circle*{4}}
 \put(60,0){\circle*{4}}
 \put(80,0){\circle*{4}}
 \put(60,20){\circle*{4}}
 \put(80,20){\circle*{4}}
 \put(0,0){\line(1,0){80}}
 \put(0,0){\line(0,1){20}}
 \put(0,20){\line(1,-1){20}}
 \put(40,0){\line(1,1){20}}
 \put(60,20){\line(1,0){20}}
 \put(80,0){\line(0,1){20}}
 \put(60,20){\line(0,-1){2}} \put(60,16){\line(0,-1){2}}
 \put(60,12){\line(0,-1){2}} \put(60,8){\line(0,-1){2}}
 \put(60,4){\line(0,-1){2}}
\end{picture}
&
\begin{picture}(80,60)(0,-20)
 \put(60,15){\circle*{4}}
 \put(80,15){\circle*{4}}
 \put(60,35){\circle*{4}}
 \put(80,35){\circle*{4}}
 \put(60,15){\line(1,0){20}}
 \put(60,15){\line(0,1){20}}
 \put(60,35){\line(1,0){20}}
 \put(80,35){\line(0,-1){20}}
 \put(0,-15){\circle*{4}}
 \put(0,5){\circle*{4}}
 \put(20,-15){\circle*{4}}
 \put(40,-15){\circle*{4}}
 \put(60,-15){\circle*{4}}
 \put(60,5){\circle*{4}}
 \put(0,-15){\line(1,0){60}}
 \put(0,-15){\line(0,1){20}}
 \put(0,5){\line(1,-1){20}}
 \put(40,-15){\line(1,1){20}}
 \put(60,5){\line(0,-1){20}}
\end{picture}
\end{tabular}
&
\begin{tabular}{cc}
\begin{picture}(80,40)(0,-20)
 \put(0,0){\circle*{4}}
 \put(0,20){\circle*{4}}
 \put(20,0){\circle*{4}}
 \put(40,0){\circle*{4}}
 \put(60,0){\circle*{4}}
 \put(40,20){\circle*{4}}
 \put(60,20){\circle*{4}}
 \put(80,20){\circle*{4}}
 \put(0,0){\line(1,0){60}}
 \put(0,0){\line(0,1){20}}
 \put(0,20){\line(1,-1){20}}
 \put(40,0){\line(0,1){20}}
 \put(40,20){\line(1,0){40}}
 \put(60,0){\line(1,1){20}}
 \put(60,20){\line(0,-1){2}} \put(60,16){\line(0,-1){2}}
 \put(60,12){\line(0,-1){2}} \put(60,8){\line(0,-1){2}}
 \put(60,4){\line(0,-1){2}}
\end{picture}
&
\begin{picture}(80,60)(0,-20)
 \put(40,15){\circle*{4}}
 \put(60,15){\circle*{4}}
 \put(40,35){\circle*{4}}
 \put(60,35){\circle*{4}}
 \put(40,15){\line(1,0){20}}
 \put(40,15){\line(0,1){20}}
 \put(40,35){\line(1,0){20}}
 \put(60,35){\line(0,-1){20}}
 \put(0,-15){\circle*{4}}
 \put(0,5){\circle*{4}}
 \put(20,-15){\circle*{4}}
 \put(40,-15){\circle*{4}}
 \put(60,-15){\circle*{4}}
 \put(60,5){\circle*{4}}
 \put(80,5){\circle*{4}}
 \put(0,-15){\line(1,0){60}}
 \put(0,-15){\line(0,1){20}}
 \put(0,5){\line(1,-1){20}}
 \put(60,5){\line(1,0){20}}
 \put(60,-15){\line(1,1){20}}
 \put(60,5){\line(0,-1){20}}
\end{picture}
\end{tabular}
\\
 {\small (4.a)
 Non induced bow tie disassembled}
&
 {\small (4.b)
 Non induced bow tie disassembled}
\\
 {\small into 1 bow tie and 1 even cycle.}
&
 {\small into 1 bow tie and 1 even cycle.}
\end{tabular}
\end{center}
 \caption{Non induced bow ties are decomposable}
 \label{figure_disassembling_bowtie}
\end{figure}

\demo Using Proposition~\ref{all_indecomposable_even_walks} and
observing that induced bow ties are indecomposable, we will be done
once is shown that every non induced bow tie is decomposable. Given
a non induced bow tie, there is at least one edge in ${\cal G}(\ff)$
which is not an edge of the bow tie and that connects two vertices
of the bow tie. Depending on the kind of vertices connected by this
extra edge, one gets four distinct situations:
\begin{itemize}
\item[(1)] one vertex on one structural odd cycle and the second on the other;
\item[(2)] one vertex on one structural odd cycle and the other on the structural path;
\item[(3)] both vertices on the structural path;
\item[(4)] both vertices on the same structural odd cycle.
\end{itemize}
Figure~\ref{figure_disassembling_bowtie} illustrates these four
situations. The dot edge is, in each situation, the extra edge that
makes the bow tie non induced. Observe that the decomposition may
depend on the parity of the number edges in some specific part of
the configuration. For example, in situation (3), depending on the
parity of the number of edges on the structural path that connect
the two vertices joined by the dot edge, one gets (3a) or (3b). In
(4), the dot edge is a chord of one of the structural odd cycles and
hence it divides it into an even cycle and an odd cycle. When the
odd cycle is connected to the structural path of the non induced bow
tie, one has (4a), otherwise one has (4b).

 \noindent
Note that in each situation, the dot edge is used as decomposing set
except in (4.b) where the decomposing set contains the dot edge and
some edges of the non induced bow tie.
 \qed

\begin{Remark}\rm Induced even cycles are certainly indecomposable
but also cycles having a chord may be indecomposable. Of course, the
existence of a chord subdividing the induced subgraph associated to
the cycle vertices into smaller even cycles makes it decomposable as
the first example in Figure~\ref{figure_decomposable} shows. But a
cycle can also be decomposable if this condition is not fulfilled as
the third example in Figure~\ref{figure_decomposable} illustrates.
\end{Remark}

\section{Combinatorics and polar syzygies}

In this section we establish the nature of generators of both the
differential syzygy module $\fz$ and its counterpart, the polar
syzygy module $\fp$ -- see Section~\ref{terminology} for the needed
terminology.

\subsection{Even closed walks induce syzygies}\label{central_results}

Recall that, as in Section~\ref{terminology}, the elements of $\fz$
(respectively, of $\fp$) are named differential (respectively,
polar) syzygies of $\ff$.

We have the following basic result.

\begin{Lemma}\label{walk_is_syzygy}
Let $\ff\subset R$ be a set of monomials of degree $2$ and let
${\mathfrak w}=\{g_1,\ldots,g_{2r}\}$ be an even closed walk of
${\cal G}(\ff)$ {\rm(}$r\geq 2${\rm)}. Then the transpose of the
vector
$$\zwt:=\left(\frac{g}{g_1},-\frac{g}{g_2}, \frac{g}{g_3},\ldots,
-\frac{g}{g_{2r}}\right)$$ is a differential syzygy of the edge
sequence $\{g_1,\ldots,g_{2r}\}$, where $g$ stands for the least
common multiple of the distinct monomials in the sequence
$g_1,\ldots,g_{2r}$.
\end{Lemma}
\demo Assume that $g_1=x_{i_1}x_{i_2}$, $g_2=x_{i_2}x_{i_3}$,
\ldots, $g_{2r}=x_{i_{2r}}x_{i_1}$. One has
\begin{eqnarray*}
\frac{g}{g_1}\,dg_1-\frac{g}{g_2}\,dg_2&=&\frac{g}{x_{i_1}}dx_{i_1}+\frac{g}{x_{i_2}}dx_{i_2}-
\left(\frac{g}{x_{i_2}}dx_{i_2}+\frac{g}{x_{i_3}}dx_{i_3}\right)\\
&=&\frac{g}{x_{i_1}}dx_{i_1}-\frac{g}{x_{i_3}}dx_{i_3}
\end{eqnarray*}
as elements of $\sum_{i=1}^n Rdx_i$. Inducting, one gets at the
$(2r-2)$nd step
$$\frac{g}{g_1}dg_1-\frac{g}{g_2}dg_2+\cdots
-\frac{g}{g_{2r-2}}dg_{2r-2}=
\frac{g}{x_{i_1}}dx_{i_1}-\frac{g}{x_{i_{2r-1}}}dx_{i_{2r-1}}.$$
Applying two more steps and recalling that
$g_{2r}=x_{i_{2r}}x_{i_1}$, it is clear that
$$\frac{g}{g_1}dg_1-\frac{g}{g_2}dg_2+\cdots-\frac{g}{g_{2r}}dg_{2r}=0.$$
\qed

\bigskip

We associate to an even closed walk ${\mathfrak
w}=\{g_1,\ldots,g_{2r}\}$ of ${\cal G}(\ff)$, a vector $\zw$ in
$R^n$ as follows: denoting by $(\zwt)_j$ the $j$-th entry of the
vector $\zwt$ defined in Lemma~\ref{walk_is_syzygy} ($1\leq j\leq
2r$), the $i$-th entry of $\zw$ ($1\leq i\leq n$) is
$\displaystyle{\sum_{j\,/\,g_j=f_i}(\zwt)_j}$ (understanding that
this is $0$ if $f_i$ does not belong to the edge sequence of
${\mathfrak w}$). Note that if the even closed walk ${\mathfrak w}$
is non-split, the $i$-th entry of $\zw$ is $0$ if and only if $f_i$
does not belong to the edge sequence of ${\mathfrak w}$ by
Lemma~\ref{edge_rep_in_non-split_walk} (3). Moreover, by
Lemma~\ref{edge_rep_in_non-split_walk} (2), the nonzero entries of
$\zw$ are pure monomials in $R$ with a factor $\pm 1$ or $\pm 2$.

\begin{Example}\label{example_villa_diffsyz}{\rm
Consider $\ff=\{f_1,\ldots,f_5\}\subset R=K[x_1,x_2,x_3]$ with
$f_1=x_1^2$, $f_2=x_1x_2$, $f_3=x_2^2$, $f_4=x_2x_3$, $f_5=x_3^2$
whose associated graph ${\cal G}(\ff)$ is shown in
Figure~\ref{figure_indecomposable}. If ${\mathfrak w}$ is the
induced looped bow tie in ${\cal G}(\ff)$ involving the first and
the third loops, then
$\zw=(x_2x_3^2,-2x_1x_3^2,0,2x_1^2x_3,-x_1^2x_2)^t\in R^{5}$.
}\end{Example}

\medskip
The following result is one of the basic bridging devices between
combinatorics and polarizability. Keeping the just introduced
notation, one has:

\begin{Theorem}\label{syzygies_are_walks}
Let $\ff\subset R$ be a set of monomials of degree $2$. Then the
differential syzygy module $\fz$ of $\ff$ is generated by the
vectors $\zw$, for all non-split even closed walks ${\mathfrak w}$
of length $\geq 4$ of the graph ${\cal G}(\ff)$.
\end{Theorem}
\demo By Lemma~\ref{walk_is_syzygy}, for every even closed walk
${\mathfrak w}=\{g_1,\ldots,g_{2r}\}$ the transpose of $\zwt$ is a
syzygy of the differentials of the edge sequence
$\{g_1,\ldots,g_{2r}\}$. Suppose that ${\mathfrak w}$ is non-split.
For any edge repetition $g_j=g_l$ in the edge sequence,  identify
the corresponding differentials $dg_j, dg_l$ and, accordingly,
introduce a factor of $\pm 2$ as coefficient of the corresponding
coordinate of $\zwt$ because $j\equiv l\,(\hbox{mod}\, 2)$. Next,
complete the transpose of $\zwt$ to a full vector of $R^m$ by
placing $0$ at every coordinate corresponding to an $f_j\notin
\{g_1,\ldots,g_{2r}\}$. In this way, the resulting vector of $R^m$
clearly belongs to $\fz$.

Conversely, let $\mbox{\sc z} \in \fz$ be a differential syzygy of
$\ff$. Since $\ff$ is a set of monomials of the same degree, the
transposed Jacobian module  $\fd(\ff)$ in its natural embedding in
$\sum_{i=1}^nR\, dx_i$ is graded with respect to the fine grading.
Therefore, it has a minimal ${\mathbb Z}^{n}$-graded free resolution
and, in particular, $\mbox{\sc z}$ is an $R$-linear combination of
vectors $\mbox{\sc z}_1,\ldots, \mbox{\sc z}_t$ in $\fz\subset R^m$
whose coordinates are terms $\alpha\xx^{\bf a}\in R$ with $\alpha\in
{\mathbb Q}$. Multiplying each $\mbox{\sc z}_i$ by an integer, one
can assume without loss of generality that any differential syzygy
is an $R$-linear combination of vectors in $\fz\subset R^m$ whose
coordinates are terms $\alpha\xx^{\bf a}\in R$ with $\alpha\in
{\mathbb Z}$. Thus, assume that the given differential syzygy
$\mbox{\sc z}$ is already of the latter form, so that one has a
relation of the form
 $
\alpha_1 \xx^{{\bf a}_1}df_1+\alpha_2 \xx^{{\bf
a}_2}df_2+\cdots+\alpha_m \xx^{{\bf a}_m}df_m\,=\,0
 $ with $\alpha_i\in {\mathbb Z}$. In other words, one can assume
 that the given differential syzygy
$\mbox{\sc z}$ gives a relation of the form
\begin{equation}\label{relation_from_syzygy}
\epsilon_1 M_1dg_1+\epsilon_2 M_2dg_2+\cdots+\epsilon_s
M_sdg_s\,=\,0
\end{equation}
where $g_1,\ldots,g_s\in \ff$,
$\epsilon_1,\ldots,\epsilon_s\in\{-1,+1\}$, and $M_1,\ldots,M_s$ are
monomials in $R$ such that $\gcd(M_1,\ldots,M_s)=1$ and $M_i=M_j$
(and $\epsilon_i=\epsilon_j$) whenever $g_{i}=g_j$ for some $1\leq
i<j\leq s$. Moreover, one can also assume that this relation is
shortest for $dg_1,\ldots,dg_s$. In this situation we claim that
$\mbox{\sc z}=\zw$ for some non-split even closed walk ${\mathfrak
w}$.

Indeed, write $g_1=x_{i_1}x_{i_2}$. Then, by the same token as in
the proof of Lemma~\ref{walk_is_syzygy}, one has
$M_1dg_1=M_1x_{i_2}\,dx_{i_1}+M_1x_{i_1}\,dx_{i_2}$  (including the
collapsing case $i_1=i_2$, whereby $M_1dg_1$ has one single non-zero
coordinate, namely, $2M_1x_{i_1}$ as coefficient of $dx_{i_1}$). Now
(\ref{relation_from_syzygy}), forces the existence of an index
$\ell$, $2\leq\ell\leq s$, such that $\epsilon_\ell=-1$ and that one
of the two non-zero coordinates of the vector $M_\ell dg_\ell$ is
$M_1x_{i_1}$ as coefficient of $dx_{i_2}$. Moreover, the other
non-zero coordinate cannot be a coefficient of $dx_{i_1}$. In other
words, upon reordering the $g_j$'s if necessary, one can assume that
$\epsilon_2=-1$, that $g_2=x_{i_2}x_{i_3}$ for some $i_3\neq i_1$,
and $M_1x_{i_1}=M_2x_{i_3}$. Then $M_1dg_1-
M_2dg_2=M_1x_{i_2}\,dx_{i_1}-M_2x_{i_2}\,dx_{i_3}$, with $i_3\neq
i_1$. By the same argument, there exists $\ell$, $3\leq\ell\leq s$,
such that $\epsilon_\ell=+1$ and with the property that one of the
non-zero coordinates of the vector $M_\ell dg_\ell$ is $M_2x_{i_2}$
as coefficient of $dx_{i_3}$. Again, upon reordering the $g_j$'s if
necessary, one can assume that $\ell=3$, i.e., $g_3=x_{i_3}x_{i_4}$
for some $i_4\neq i_2$, and $M_2x_{i_2}=M_3x_{i_4}$. Then $M_1dg_1-
M_2dg_2+M_3dg_3=M_1x_{i_2}\,dx_{i_1}+M_3x_{i_3}\,dx_{i_4}$.
Iterating this process and reordering the $g_i$'s at each step if
necessary, one gets that $g_j=x_{i_j}x_{i_{j+1}}$ for all
$j=1,\ldots,s$. Note that in order to get the zero vector, $s$ has
to be even, and $i_{j+1}=i_1$. In other words, $\{g_1,\ldots,g_s\}$
is an even closed walk. Moreover, the condition that has to be
satisfied by the monomials $M_j$ at each step is
\begin{equation}\label{cond_on_monomials_in_diff_syzygy}
M_j\frac{g_j}{\gcd(g_j,g_{j+1})}=
M_{j+1}\frac{g_{j+1}}{\gcd(g_j,g_{j+1})},\ \forall\,j=1,\ldots,s\ .
\end{equation}
Setting $M:=M_1 g_1$, one has that $M=M_j g_j$ for all
$j=1,\ldots,s$. Moreover,
 $\displaystyle{\frac{g_{j+1}}{\gcd(g_j,g_{j+1})}}$
divides $M_j$, and hence $\lcm(g_j,g_{j+1})$ divides $M$. Letting
$g$ stand for the least common multiple of the distinct monomials in
the sequence $g_1,\ldots,g_s$, this implies the existence of a
monomial $N\in R$ such that $M=g N$. Then, for all $j=1,\ldots,s$,
$\displaystyle{M_j=\frac{g}{g_j}N}$. Since we have assumed that the
monomials $M_j$ have no non-trivial common factor, one has that
$N=1$, and hence $\mbox{\sc z}=\zw$ for the even closed walk
${\mathfrak w}:=\{g_1,\ldots,g_s\}$.

Finally, observe that if an even closed walk ${\mathfrak w}$ splits
into two smaller even closed walks ${\mathfrak w}_1$ and ${\mathfrak
w}_2$, then
 $\displaystyle{
 \zw=\frac{g}{\ell_1}\mbox{\sc z}_{{\mathfrak w}_1}+
     \frac{g}{\ell_2}\mbox{\sc z}_{{\mathfrak w}_2}
 }$
where $g$, $\ell_1$ and $\ell_2$ are the least common multiples of
the monomials in the edge sequences associated to ${\mathfrak w}$,
${\mathfrak w}_1$ and ${\mathfrak w}_2$ respectively. \qed

\bigskip

Let $P\subset k[\TT]$ be the presentation ideal of $k[\ff]$ relative
to the given generators $\ff$. We formally introduce a construct
that is a special polar syzygy to play a central role in the
discussion.

\begin{Definition}\label{defining_tw}\rm
Let ${\mathfrak w}$ denote an even closed walk  of the graph ${\cal
G}(\ff)$.  To it one associates the binomial relation $p_{\mathfrak
w}=\TT_{{\mathfrak w}^+}-\TT_{{\mathfrak w}^-}\in P$ in a notation
mimicking that of \cite[7.1.4]{VillaBook}. Define the {\em
associated polar syzygy} $\tw$ to be the differential of
$p_{\mathfrak w}$ further evaluated at the edges of ${\mathfrak w}$.
In further detail,  regarding $\tw$ as a column vector, its $j$th
coordinate is the $T_j$th derivative of $p_{\mathfrak w}$ (hence, a
monomial) further evaluated at the corresponding edge $g_j$ in the
edge sequence of the walk ${\mathfrak w}$.
\end{Definition}

\begin{Example}\label{example_villa_polarsyz}{\rm
If ${\mathfrak w}$ is the induced bow tie considered in
Example~\ref{example_villa_diffsyz}, then $p_{\mathfrak
w}=T_1T_4^2-T_2^2T_5$. The associated polar syzygy is
$\tw=(x_2^2x_3^2,-2x_1x_2x_3^2,0,2x_1^2x_2x_3,-x_1^2x_2^2)^t$. Note
that this polar syzygy is related to the differential syzygy $\zw$
determined in Example~\ref{example_villa_diffsyz} by $\tw=x_2\zw$.
As we shall argue in Lemma~\ref{tw_vs_zw}, this relation is not
accidental. }\end{Example}

Of a similar nature is the following counterpart to
Theorem~\ref{syzygies_are_walks}.

\begin{Theorem}\label{polarsyzygies_are_walks}
Let $\ff\subset R$ be a set of monomials of degree $2$. Then the
polar syzygy module $\fp$ of $\ff$ is generated by the vectors
$\tw$, for all non-split even closed walks ${\mathfrak w}$ of length
$\geq 4$ of the graph ${\cal G}(\ff)$.
\end{Theorem}

\demo Since there is no particular claim about a minimal set of
generators, it will suffice to argue that: (1) the presentation
ideal $P$ as above is generated by the polynomials $p_{\mathfrak
w}$, for all even closed walks ${\mathfrak w}$ of the graph ${\cal
G}(\ff)$; (2)  if an even closed walk ${\mathfrak w}$ splits into
smaller cycles ${\mathfrak w}_1$ and ${\mathfrak w}_2$, then the
corresponding polynomial $p_{\mathfrak w}$ is superfluous in the
sense that it belongs to the subideal generated by the polynomials
$p_{{\mathfrak w}_1}$ and $p_{{\mathfrak w}_2}$.

We deal with the second claim first as it is visible offhand.
Namely, one has in the previous notation $ p_{\mathfrak w}
 =(\TT_{{\mathfrak w}_2^+})(\TT_{{\mathfrak w}_1^+}-\TT_{{\mathfrak w}_1^-})
 +(\TT_{{\mathfrak w_1}^-})(\TT_{{\mathfrak w}_2^+}-\TT_{{\mathfrak w}_2^-})
 =(\TT_{{\mathfrak w}_2^+})p_{{\mathfrak w}_1}
 +(\TT_{{\mathfrak w_1}^-})p_{{\mathfrak w}_2}
$.

As for the first claim, we note that it is \cite[Proposition~8.1.2
(a)]{VillaBook} when the graph ${\cal G}(\ff)$ is simple. In
general, if loops are taken into consideration, the same proof works
with minor adaptation. Indeed, setting ${\cal B}:=\{p_{\mathfrak
w}\,|\, {\mathfrak w}\hbox{ is an even closed walk}\}$, one has
$({\cal B})\subset P$ as already pointed out before. Denoting by
$P_s$ the
 part of the toric ideal $P$ of degree $s$, we show by induction on $s\geq 2$
 that $P_s\subset ({\cal B})$. Thus, let $p\in P_2$ be any binomial, say
$p=T_{i_1}T_{i_2}-T_{i_3}T_{i_4}$ for $1\leq i_1,i_2,i_3,i_4\leq m$
with $i_1\neq i_3,i_4$ and $i_2\neq i_3,i_4$. At least one of the
monomials $f_{i_1}, f_{i_2}, f_{i_3}, f_{i_4}$ is square free
(otherwise $f_{i_1}=f_{i_3}$ or $f_{i_1}=f_{i_4}$). One can assume
without loss of generality that $f_{i_1}=x_1x_2$, that $x_1$ divides
$f_{i_3}$ and that $x_2$ divides $f_{i_4}$, i.e., $f_{i_3}=x_1x_j$
and $f_{i_4}=x_2x_k$ for $1\leq j,k\leq n$ such that $j\neq 2$ and
$k\neq 1$. If $j=1$ and $k=2$, then $p=p_{\mathfrak w}$ where
 ${\mathfrak w}$ is a monedge bow tie whose structural cycles are loops.
If $j=1$ (and $k\neq 2$), or $k=2$ (and $j\neq 1$), or $j=k$ (and
$j\neq 1$, $k\neq 2$), then $p=p_{\mathfrak w}$ where
 ${\mathfrak w}$ is a path-degenerate looped bow tie whose structural
cycle is a $3$-cycle. Finally, if $j\neq 1$, $k\neq 2$ and $j\neq
k$, $p=p_{\mathfrak w}$ where ${\mathfrak w}$ is a $4$-cycle. Thus,
$P_2\subset ({\cal B})$. In order to show that $P_s\subset ({\cal
B})$ once we assume that $P_t\subset ({\cal B})$ for all $t<s$, we
use an argument similar to the one used in loc. cit. when the graph
${\cal G}(\ff)$ is simple. Let $p=T_{i_1}\cdots
T_{i_s}-T_{j_1}\cdots T_{j_s}$ be a binomial in $P_s$ with $1\leq
i_1,\ldots,i_s,j_1,\ldots,j_s \leq m$. If, relabeling the
generators, one has that $f_{i_1}\cdots f_{i_r}=f_{j_1}\cdots
f_{j_r}$ for some $r<s$, then the relation $p=T_{i_{r+1}}\cdots
T_{i_s}(T_{i_1}\cdots T_{i_r}-T_{j_1}\cdots T_{j_r})+T_{j_1}\cdots
T_{j_r}(T_{i_{r+1}}\cdots T_{i_s}-T_{j_{r+1}}\cdots T_{j_s})$ and
the induction hypothesis imply that $p\in ({\cal B})$. Assume now
that $f_{i_1}\cdots f_{i_r}\neq f_{j_1}\cdots f_{j_r}$ for all $r<s$
and any relabeling of the elements $f_{i_1}, \ldots, f_{i_s},
f_{j_1}, \ldots, f_{j_s}$. Since $f_{i_1}\cdots f_{i_s}=
f_{j_1}\cdots f_{j_s}$, relabeling $f_{j_1},\ldots f_{j_s}$ is
necessary, one can assume without loss of generality that
$f_{i_1}=x_{k_1}x_{l_1}$ and $f_{j_1}=x_{l_1}x_{k_2}$ for $1\leq
k_1,l_1,k_2\leq n$ such that $k_1\neq k_2$ (note that if $f_{i_1}$,
respectively $f_{j_1}$, corresponds to a loop in ${\cal G}(\ff)$,
then $k_1=l_1$, respectively $l_1=k_2$). Thus, $x_{k_2}$ divides
$f_{i_2}\cdots f_{i_s}$ and one can assume that
$f_{i_2}=x_{k_2}x_{l_2}$ for some $1\leq l_2\leq n$. One has that
$f_{i_1}f_{i_2}=x_{k_1}x_{l_2}f_{j_1}$, and hence $x_{l_2}$ divides
$f_{j_2}\cdots f_{j_s}$. At this step, one has that
$f_{i_1}=x_{k_1}x_{l_1}$, $f_{j_1}=x_{l_1}x_{k_2}$,
$f_{i_2}=x_{k_2}x_{l_2}$, and one can assume that
$f_{j_2}=x_{l_2}x_{k_3}$ for $k_3\neq k_1$ unless $s=2$. Iterating
the argument, one gets an even closed walk ${\mathfrak w}$ of the
graph ${\cal G}(\ff)$ such that $p=p_{\mathfrak w}$. \qed

\bigskip

Next we clarify the precise relation between the polar syzygy $\tw$
and its differential counterpart $\zw$, for a given non-split even
closed walk ${\mathfrak w}$.

\begin{Lemma}\label{tw_vs_zw}
Let ${\mathfrak w}$ denote a non-split even closed walk on a graph
${\cal G}(\ff)$. Then,
$$\tw=M\zw\,,$$
where $M$ is the product of the repeated vertices in the closed walk
obtained by removing from ${\mathfrak w}$ the loops {\rm(}$M=1$ if
it has no vertex repetition{\rm)}. In particular, this applies to
the following particular configurations:
\begin{enumerate}
\item[{\rm (1)}]
If ${\mathfrak w}$ is either an even cycle, a path-degenerate looped
bow tie, or a monedge bow tie whose structural cycles are loops,
then $\tw=\zw$.
\item[{\rm (2)}]
If ${\mathfrak w}$ is a path-degenerate bow tie which is not looped,
and if $x_i$ is the common vertex of its two structural cycles, then
$\tw=x_i\zw$.
\item[{\rm (3)}]
If ${\mathfrak w}$ is a bow tie which is neither path-degenerate nor
a monedge bow tie whose structural cycles are loops, and if $N$ is
the product of the vertices of the structural connecting path
excluding the base vertex of the structural odd cycle when the later
is a loop, then $\tw=N\zw$.
\end{enumerate}
\end{Lemma}

\demo Consider a non-split even closed walk ${\mathfrak
w}=\{g_1,\ldots,g_{2r}\}$ on ${\cal G}(\ff)$. On one hand, recall
that in order to get $\tw$ one takes $\TT$-derivatives of the
binomial $p_{{\mathfrak w}}=T_1T_3\cdots T_{2r-1}-T_2T_4\cdots
T_{2r}$ and evaluate every $T_j$ on the corresponding edge $g_j$ in
the edge sequence of the walk (see Definition~\ref{defining_tw}) --
as a slight check, note that the $\TT$-degree of $p_{\mathfrak w}$
is $r$, hence the $\xx$-degree of $\tw$ is the even integer
$2(r-1)=2r-2$. Thus, typically, the first coordinate reads
$g_3g_5\cdots g_{2r-1}$ (respectively,
$2g_1g_3g_5\cdots\widehat{g_{2j+1}}\cdots g_{2r-1}= 2g_3g_5\cdots
g_{2r-1}$) if $g_1$ is not repeated (respectively, if $g_1=g_{2j+1}$
for some $j\geq 1$). On the other hand, the least common multiple of
$\{g_1,\ldots,g_{2r}\}$ is the monomial
$$
g
 =\frac{x_{i_1}\cdots x_{i_{2r}}}{M}=\frac{g_1g_3\cdots g_{2r-1}}{M}
 =\frac{g_2g_4\cdots g_{2r}}{M}
$$
and one readily obtains the required relation.

Of course, (1), (2) and (3) follow readily from the general
statement.
 \qed

\subsection{Minimal sets of generators}

In this part we seek to squeeze down the previous slightly loose
sets of generators to minimal sets of generators of the modules
$\fz$ and $\fp$.

First a word about sets of minimal generators of these modules.
Since $\fz$ is the module of syzygies of the transposed Jacobian
module $\fd(\ff)\subset \sum_{i=1}^nR\, dx_i$ and the latter is a
graded $k[\xx]$-module with respect to the standard graded
polynomial ring $k[\xx]$, then $\fz$ is a graded submodule, say,
$\fz=\oplus_{s^\geq 0} \fz_s$. Theorem~\ref{syzygies_are_walks}
tells us a set of generators of $\fz$. This set can in theory be
squeezed to a minimal set of generators; for any such set of
generators, the number of elements in each graded piece $\fz_s$ is
invariant as is the highest possible degree $s$ of an element in it.
We will agree to say that an element of $\fz$ is {\em superfluous}
in the sense that it does not belong to any set of minimal
generators of $\fz$ (not just lying outside a specific such set).
This is of course the counterpart to the usual notion of an absolute
minimal generator $z$ which, in our setting, reads as $z\in
\fz\setminus (\xx)\fz$. Of course, a test for knowing that $z$ is
superfluous is that its degree be larger than the uniquely defined
highest generation degree of the module $\fz$; however, we in
general have no theoretic hold of this degree.

A similar phenomenon happens in $\fp$ as the latter is in its turn a
graded submodule of $\fz$. Here one can pretest superfluity of an
element $\tw$ in $\fp$ by testing whether the associated binomial
$p_{\mathfrak w}$ is a minimal generator of the defining ideal $P$.
Unfortunately this works only in one direction in general (see
Remark~\ref{mu_fp_mu_P}).

We give some examples to illustrate this order of ideas in our
present setting, stressing additionally that even non-split or
indecomposable walks may be (absolute) superfluous.

\begin{Example}\label{minimality_in_fz}\rm Consider the simple
graph in Figure~\ref{figure_hexagone}. Here the hexagon ${\mathfrak
w}$ is non-split, but the corresponding $\zw$ is deep inside the
submodule generated by the vectors corresponding to the square and
the path-degenerate bow tie and these two form a set of minimal
generators of $\fz$. How do we know that $\zw$ is (absolute)
superfluous? Simply because, by definition, its degree is $4$ which
is larger than the generation degree $3$ of $\fz$. As an additional
remark, since $\zw=\tw$ for an hexagon and $\tw$ is part of a
minimal set of generators of $\fp$, then the edges of the graph do
not form a polarizable set.

\begin{figure}[h]
\setlength{\unitlength}{.032cm}
\begin{center}

\begin{picture}(60,30)(0,10)
 \put(0,0){\circle*{4}}
 \put(0,35){\circle*{4}}
 \put(30,45){\circle*{4}}
 \put(30,-10){\circle*{4}}
 \put(60,35){\circle*{4}}
 \put(60,0){\circle*{4}}
 \put(0,0){\line(0,1){35}}
 \put(30,-10){\line(-3,1){30}}
 \put(30,-10){\line(3,1){30}}
 \put(30,45){\line(-3,-1){30}}
 \put(30,45){\line(3,-1){30}}
 \put(60,0){\line(0,1){35}}
 \put(60,35){\line(-1,0){60}}
 \put(30,-10){\line(2,3){30}}
\end{picture}
\end{center}
 \caption{Non-split even closed walk providing a superfluous differential syzygy}
 \label{figure_hexagone}
\end{figure}
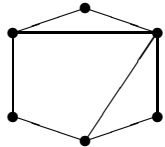
\end{Example}

\begin{Example}\label{minimality_in_fp}\rm For an example in $\fp$,
consider the graph in Figure~\ref{figure_indecomposable}. The even
closed walk ${\mathfrak w}$ of length $6$ supported by the bow tie
involving the first and third loops is such that $p_{\mathfrak w}$
belongs to the ideal $(p_{{\mathfrak w}_1},p_{{\mathfrak
w}_2})\subset k[\TT]$, where ${\mathfrak w}_1, {\mathfrak w}_2$ are
the two even closed walks of length $4$ supported by the other two
bow ties. By a previous observation above $\tw$ is not an absolute
minimal generator of $\fp$. Note that ${\mathfrak w}$ is non-split
(and even indecomposable as observed right after
Figure~\ref{figure_indecomposable}), and that $\zw$ is a actually
minimal generator of $\fz$; of course, necessarily, $\tw\neq \zw$.
\end{Example}

\medskip

One can now improve on the result of
Theorem~\ref{syzygies_are_walks} as a first approximation to
describing a minimal generating set of the differential syzygy
module ${\fz}$.

\begin{Theorem}\label{minimal_gen_set_of_syzygies}
Keeping the previous notation, the syzygy module ${\fz}$ of ${\cal
D}(\ff)$ is generated by the vectors $\zw$, for all even cycles and
induced bow ties ${\mathfrak w}$ of the graph ${\cal G}(\ff)$.
\end{Theorem}
\demo By Theorem~\ref{syzygies_are_walks}, ${\fz}$ is generated by
the vectors $\zw$, for all non-split even closed walks ${\mathfrak
w}$ of the graph ${\cal G}(\ff)$ of length $\geq 4$.

We first show that if ${\mathfrak w}$ is a non-split even closed
walk that contains a cycle of which at least two vertices are vertex
repetitions of ${\mathfrak w}$, then $\zw\in (\xx)\fz$. This
statement is proved using the decomposition used in the proof
Lemma~\ref{lemma_on_number_of_cycles_in_indecomposable}: in this
situation one can readily check that if ${\mathfrak w}_1$ and
${\mathfrak w}_2$ are the even closed walks introduced there, then
 $\displaystyle{
 \zw=\frac{g}{\ell_1}\mbox{\sc z}_{{\mathfrak w}_1}+
     \frac{g}{\ell_2}\mbox{\sc z}_{{\mathfrak w}_2}
 }$
where $g$, $\ell_1$ and $\ell_2$ are the least common multiples of
the monomials in the edge sequences associated to ${\mathfrak w}$,
${\mathfrak w}_1$ and ${\mathfrak w}_2$ respectively. Since the sets
of variables involved in the vertex sequences of ${\mathfrak w}_1$
and ${\mathfrak w}_2$ are both strictly contained in the set of
variables involved in the vertex sequence of ${\mathfrak w}$, one
has that
 $\displaystyle{
 \frac{g}{\ell_1}, \frac{g}{\ell_2}
 \neq 1
 }$,
and hence $\zw\in (\xx)\fz$.

As a consequence, following the argument in the proof of
Proposition~\ref{all_indecomposable_even_walks} we deduce that
${\fz}$ is at least generated by the vectors $\zw$, for all even
cycles and bow ties ${\mathfrak w}$ of the graph ${\cal G}(\ff)$.

To complete the proof we  show that this set of generators can be
further shrunk. Namely, we now show that if ${\mathfrak w}'$ is a
non induced bow tie, then the differential syzygy $\mbox{\sc
z}_{{\mathfrak w}'}$ belongs to the submodule generated by the
vectors $\zw$ for all cycles and induced bow ties ${\mathfrak w}$ of
${\cal G}(\ff)$. We induct on the number of the induced edges of the
graph adjacent to vertices of ${\mathfrak w}'$,  off the structural
edges of ${\mathfrak w}'$. If this number is zero - i.e., no
additional such edges, then the bow tie is non induced, hence the
result is vacuously satisfied.

In order to apply the inductive hypothesis, refer back to the
decomposition of ${\mathfrak w}'$ into two even closed walks
${\mathfrak w}_1$ and ${\mathfrak w}_2$ as in the proof of
Proposition~\ref{indecomposable_are_cycles_and_induced_bowties}.
Note that this provides a relation
 $\displaystyle{
 \mbox{\sc z}_{{\mathfrak w}'}
 =\frac{g}{\ell_1}\mbox{\sc z}_{{\mathfrak w}_1}+
     \frac{g}{\ell_2}\mbox{\sc z}_{{\mathfrak w}_2}
 }$
where $g$, $\ell_1$ and $\ell_2$ are the least common multiples of
the monomials along the structural edge sequences of ${\mathfrak
w}'$, ${\mathfrak w}_1$ and ${\mathfrak w}_2$ respectively. This
holds for any of the basic ways described in
Figure~\ref{figure_disassembling_bowtie} in which a non induced bow
tie can decompose. Now, with one single exception, ${\mathfrak w}_1$
and ${\mathfrak w}_2$ are even cycles or bow ties. The exception is
when, say, ${\mathfrak w}_1$ is an even molecule (see
Figure~\ref{figure_disassembling_bowtie}, (3.b)). But then
${\mathfrak w}_2$ is an even cycle, and the molecule ${\mathfrak
w}_1$ again decomposes further into a bow tie and an even cycle
which is ${\mathfrak w}_2$.

Thus, in all situations, one has
 $\displaystyle{
 \mbox{\sc z}_{{\mathfrak w}'}
 =\lambda\frac{g}{\ell_1}\mbox{\sc z}_{{\mathfrak w}_1}+
     \frac{g}{\ell_2}\mbox{\sc z}_{{\mathfrak w}_2}
 }$
where ${\mathfrak w}_1, {\mathfrak w}_2$ are even cycles or bow
ties, $g$, $\ell_1$ and $\ell_2$ are the least common multiples of
the monomials along the structural edge sequences of  ${\mathfrak
w}'$, ${\mathfrak w}_1$ and ${\mathfrak w}_2$ respectively, and
$\lambda=1$ except in the basic situation (3.b) where $\lambda=2$,
${\mathfrak w}_1$ is an even cycle and ${\mathfrak w}_2$ a bow tie.

If, say, ${\mathfrak w}_1$  is a non induced bow tie,  the number of
the induced edges of the graph adjacent to vertices of ${\mathfrak
w}_1$,  off the structural edges of ${\mathfrak w}_1$, is strictly
smaller than the analogous number corresponding to ${\mathfrak w}'$.
Therefore, we can apply the inductive hypothesis and the result
follows suit. \qed

\begin{Remark}\label{indecomposable_not_sufficient}\rm
In general, one cannot replace even cycles by indecomposable even
cycles in Theorem~\ref{minimal_gen_set_of_syzygies}. Consider the
graph whose edges $\ff$ are those of a decagon, i.e., a $10$-cycle
with  vertices labeled $x_1,\,x_2,\,\ldots,\,x_{10}$, and in
addition the chords $x_2x_8$ and $x_3x_7$. A straightforward
calculation shows that $\ff$ is polarizable -- see also
Theorem~\ref{polar_characterization}. Moreover, the differential
syzygy module is minimally generated by the $4$-cycle
$\{x_2,\,x_8,\,x_7,\,x_3,\,x_2\}$ and the entire $10$-cycle. To be
in conformity with the result of
Theorem~\ref{minimal_gen_set_of_syzygies}, note that the  monedge
bow tie whose structural odd cycles are both of length $5$
 and whose structural path is the edge $x_2x_3$ is
decomposable -- using as decomposing set the edge $x_7x_8$ as in
Figure~\ref{figure_disassembling_bowtie}, (1). On the other hand,
the $10$-cycle is decomposable with decomposing set the edges
$x_2x_8$, $x_2x_3$ and $x_3x_7$, by which it disassembles into the
$4$-cycle and the monedge bow tie. \end{Remark}

\begin{Remark}\label{mu_fp_mu_P}\rm A point that would require further clarification is a criterion
for the inequality $\mu(\fp)\leq \mu(P)$ to be an equality. An
example where a decomposable even closed walk provides a superfluous
generator of $\fp$ while the binomial $p_{\mathfrak w}$ is a
non-superfluous generator of the presentation ideal $P$ of $k[\ff]$
is illustrated by the graph of Example~\ref{bounding_cycle}. In this
example $\ff$ is polarizable. The cycle arrangement in
\cite[Ex.~8.4.14]{VillaBook} provides us with the same phenomenon
and is moreover non-polarizable.
\end{Remark}

The following result soups-up the previous result by capturing a
class of even closed walks ${\mathfrak w}$ whose associated syzygies
 $\tw$ (respectively, $\zw$) are part of a minimal set of generators
 of $\fp$ (respectively, $\fz$).

\begin{Lemma}\label{some_are_non-split}
If ${\mathfrak w}$ is an induced bow tie on a graph ${\cal G}(\ff)$,
then the associated syzygy $\tw$ {\rm (}respectively, $\zw${\rm )}
is part of a minimal set of generators of  ${\fp}$ {\rm
(}respectively, $\fz${\rm )}.
\end{Lemma}
\demo By Theorem~\ref{polarsyzygies_are_walks}, $\fp$ is generated
by the set of syzygies $\tw$ where  ${\mathfrak w}$ runs through the
set of non-split even closed walks. Since $\fp$ is a graded
$k[\xx]$-module, this set contains a subset $\mathcal{M}$ forming a
minimal set of generators of $\fp$. We claim that if ${\mathfrak w}$
is a non induced bow tie in  ${\cal G}(\ff)$ then ${\mathfrak w}\in
\mathcal{M}$.

One can assume without loss of generality that the first monomial
$f_1$ in $\ff$ corresponds to the first edge in the edge sequence of
${\mathfrak w}$; in particular, the first coordinate of the vector
$\tw$ is nonzero.

Suppose then that ${\mathfrak w}\not\in \mathcal{M}$. Write
${\mathfrak w}$ as a $k[\xx]$-linear combination of $\mbox{\sc
t}_{{\mathfrak w}_1}, \ldots, \mbox{\sc t}_{{\mathfrak w}_\ell}$,
where ${\mathfrak w}_1, \ldots, {\mathfrak w}_\ell\in \mathcal{M}$.
Then, $f_1$ belongs to the edge sequence of at least one of those
even closed walks, say ${\mathfrak w}_1=\{f_1,g_2,\ldots,g_{2r}\}$
with $g_2,\ldots,g_{2r}\in\ff$, and the first coordinate of
$\mbox{\sc t}_{{\mathfrak w}_1}$ divides the first coordinate  of
$\tw$.

Now, since ${\mathfrak w}_1$ does not coincide with ${\mathfrak w}$
because we are assuming that ${\mathfrak w}\notin \mathcal{M}$ and
since ${\mathfrak w}$ does not contain any proper even closed
 subwalk because it is a bow tie, it follows that at least one of
the monomials in the edge sequence of ${\mathfrak w}_1$, say $g_i$
for some $i\in\{2,\ldots,2r\}$, does not belong to the edge sequence
of ${\mathfrak w}$. If $g_i=x_jx_k\in\ff$  for $1\leq j\leq k\leq
n$, we claim that both $x_j$ and $x_k$ belong to the vertex sequence
of ${\mathfrak w}$. If $x_j$ does not belong to the vertex sequence
of ${\mathfrak w}$, $x_j$ does not divide any of the non zero
coordinates of $\tw$, in particular it does not divide its first
coordinate. On the other hand, recall that the first coordinate of
$\mbox{\sc z}_{{\mathfrak w}_1}$ is
 $\displaystyle{
 \frac{g}{f_1}
 }$
where $g$ stands for the least common multiple of
$f_1,g_2,\ldots,g_{2r}$, and hence $x_j$ divides the first
coordinate of the vector $\mbox{\sc z}_{{\mathfrak w}_1}$ (it
divides $g_i$ and does not divide $f_1$ because it does not belong
to the vertex sequence of ${\mathfrak w}$). By Lemma~\ref{tw_vs_zw},
$\mbox{\sc t}_{{\mathfrak w}_1}=M\mbox{\sc z}_{{\mathfrak w}_1}$ for
some monomial $M\in R$, and hence $x_j$ divides the first coordinate
of $\mbox{\sc t}_{{\mathfrak w}_1}$ which in turn divides the first
coordinate of the vector $\tw$, a contradiction.

We have thus shown that there exists an edge $x_jx_k\in\ff$,  not
belonging to the edge sequence of ${\mathfrak w}$, and such that
both $x_j$ and $x_k$ belong to the vertex sequence of ${\mathfrak
w}$. Therefore  the bow tie ${\mathfrak w}$ is non induced. This
wraps up the proof for a polar syzygy $\tw$.

The proof for $\zw$ is similar by drawing upon a set of generators of $\fz$
such as given in Theorem~\ref{syzygies_are_walks}. \qed

\medskip

We can now give a complete combinatorial characterization of
polarizability.

\begin{Theorem}\label{polar_characterization}
A set $\ff\subset R$ of monomials of degree $2$ is polarizable if
and only if every induced bow tie of the associated graph ${\cal
G}(\ff)$ is one of the following:
\begin{enumerate}
\item[{\rm (1)}]
A monedge bow tie whose structural cycles are loops{\rm ;}
\item[{\rm (2)}]
A path-degenerate looped bow tie.
\end{enumerate}
In particular, if $\ff$ consists only of squarefree monomials --
i.e., if the graph ${\cal G}(\ff)$ is simple -- then $\ff$ is
polarizable if and only if ${\cal G}(\ff)$ does not have any induced
bow tie.
\end{Theorem}

\demo By Lemma~\ref{tw_vs_zw}, a bow tie supporting an even closed
walk ${\mathfrak w}$ satisfies $\tw=\zw$ if and only if it is one of
the types (1) or (2) in the present statement.

Now assume that the only induced bow ties in ${\cal G}(\ff)$ are of
these types. Then any generator $\zw$ of ${\fz}$ as in
Theorem~\ref{minimal_gen_set_of_syzygies} belongs to the polar
syzygy module $\fp$, hence $\ff$ is polarizable.

Conversely, let ${\mathfrak w}$ be an even closed walk in ${\cal
G}(\ff)$ supported by an induced bow tie.  Again if ${\mathfrak w}$
is neither of the two types in the statement then $\tw=M\zw$ for
some monomial $M\neq 1$. By Lemma~\ref{some_are_non-split}, $\tw$ is
part of a minimal set of generators of $\fp$, hence $\tw\not\in
(\xx)\fp$. Therefore we must conclude that $\zw\notin \fp$, hence
$\ff$ is not polarizable. \qed

\section{Applications}

\subsection{Veronese, squarefree Veronese, bipartite}

\begin{Corollary}\label{veronese_is_polar}
Let $\ff\subset R$ be either the set of all monomials of degree $2$,
or the set of all squarefree monomials of degree $2$. Then $\ff$ is
polarizable.
\end{Corollary}

\demo In both cases, the result is a direct consequence of
Theorem~\ref{polar_characterization}. We first treat the squarefree
case. The corresponding graph is a complete simple graph (no loops).
In particular it has no induced bow ties as any induced subgraph of
a complete graph is itself complete. As for the $2$-Veronese
embedding, the corresponding graph is a complete graph with a loop
based at every vertex. This clearly forces any induced bow tie to be
either a triangle with a loop based at one of its vertices or two
loops connected by an edge. \qed

\medskip

Another consequence is a more conceptual proof of one of the main
results of \cite{jac}.

\begin{Corollary}\label{bipartite_is_polar}
Let ${\cal G}(\ff)$ denote a connected bipartite graph on edges
$\ff$. Then $\ff$ is polarizable.
\end{Corollary}
\demo It follows from Theorem~\ref{polar_characterization} because
${\cal G}(\ff)$ has no odd cycle by \cite[Prop.~6.1.1]{VillaBook}
and hence has no bow ties. \qed

\bigskip

Connected bipartite graphs admit various characterizations in
the graph literature and also in algebraic combinatorics (see, e.g.,
\cite{itg}). We will next give yet another characterization based
solely on the underlying edge-algebra. We say that a graph with
loops is connected if the underlying graph removing all loops is a
connected simple graph.

Let $A=k[\ff]\subset R$ with $\ff$ distinct monomials of degree $2$
in $n\geq 2$ variables. To any two variables $x_i,x_j$ with $i\neq
j$ we associate the $k$-algebra surjection
\begin{eqnarray*}\pi_{i,j}:R &\surjects &
S=k[x_1,\ldots,\widehat{x_j},\ldots, x_n]\\
\pi_{i,j}(x_k)&=&x_k \;\;(k\neq j)\\
\pi_{i,j}(x_j)&=&x_i
\end{eqnarray*}
Clearly, $\ker(\pi_{i,j})=(x_j-x_i)$. Set $B:=\pi_{i,j}(A)\subset S$
for the image of the restriction of this map to the $k$-subalgebra
$A$. Then $B$ is generated by the images of $\ff$, hence is still
generated by monomials $\ff'$ of degree $2$.

If ${\cal G}={\cal G}(\ff)$ and ${\cal G}'={\cal G}(\ff')$ denote
the respective associated graphs (with loops) then we say that the
corresponding graph-theoretic process is an {\it edge-pinching
operation\/} (see \cite[Corollary 4.9]{SimisVilla} where this notion
has been considered in a special case).

\begin{Proposition}\label{pinching}
Let ${\cal G}$ be a connected graph, possibly with loops, and let
$A\subset k[x_1,\ldots,x_n]$ $(n\geq 2)$ denote its associated
edge-algebra. Let ${\cal G}'$ denote the graph obtained by an
edge-pinching operation on any proper edge {\rm (i.e., not a loop)}.
Then ${\cal G}$ is bipartite {\rm(}in particular, has no loops{\rm)}
if and only if the corresponding restriction map $A \lar S$ is
injective.
\end{Proposition}

We need the following technical result.

\begin{Lemma}\label{loop_dimension}
Let ${\cal G}$ be a connected graph with $n$ vertices having at
least one loop, and let $A\subset R=k[x_1,\ldots,x_n]$ denote its
associated edge-algebra. Then $\dim A=n$.
\end{Lemma}
\demo Fix a loop, say, $x_1^2\in A$. Let $\widetilde{\cal G}$ denote
the graph obtained from ${\cal G}$ by keeping all vertices and
removing the loop $x_1^2$. Clearly, $\widetilde{\cal G}$ is still
connected. If it still contains a loop we are done by induction on
the total number of edges and loops. Thus, we may assume that
$\widetilde{\cal G}$ has no loops. If $\widetilde{\cal G}$ is not
bipartite then its associated edge-ideal has dimension $n$, hence so
does $A$. If $\widetilde{\cal G}$ is bipartite, its log-matrix $M$
has rank $n-1$. Therefore, by adding further a column
$(0,1,\ldots,0)^t$ increases the rank of $M$ by one, hence the
log-matrix of ${\cal G}$ is at least $n$, as required. \qed

\medskip

\noindent{\bf Proof of the proposition}. Suppose first that $\pi_{|A}:A \lar S$ is injective.
In particular $\pi$ does not collapse two distinct generators (edges) of $A$, hence the
images of the generators are all distinct and correspond to a graph with at least one
loop (e.g., $x_{n-1}^2)$ whose associated edge-algebra is $\pi(A)$. Clearly, this graph
is still connected. By Lemma~\ref{loop_dimension}, $\dim\pi(A)=n-1$. But then $\dim A=n-1$
as well. In particular, again by Lemma~\ref{loop_dimension}, $\widetilde{\cal G}$ has no loops,
hence must be a bipartite graph.

Conversely, if ${\cal G}$ is bipartite then $\dim A=n-1$. Once more
by Lemma~\ref{loop_dimension}, $\dim\pi(A)=n-1$. But then the
restriction $\pi_{|A}:A \surjects \pi(A)$ must have kernel zero
since $A$ is a domain. \qed

\medskip

The following result shows that polarizability depends on the chosen
embedding $A\subset R$, hence is not an invariant property of the
algebra $A$.

\begin{Proposition}\label{pinching_chordal_evengon}
Let ${\cal G}$ be a graph {\rm (or an induced subgraph)} consisting
of an even cycle with one single chord inducing a decomposition in
two smaller even cycles. Then the graph ${\cal G}'$ obtained by
pinching the chord {\rm (see Figure~\ref{figure_chordal_evengon}
below)} is not polarizable.
\end{Proposition}

\demo By edge-pinching we have created an induced path-degenerate
bow tie whose structural cycles are not loops and the result follows
from Theorem~\ref{polar_characterization}. \qed

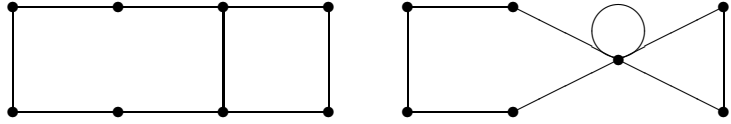
\begin{figure}[h]
\begin{center}
\vskip 0.5cm \setlength{\unitlength}{.035cm}
\begin{picture}(270,40)(0,0)
 \put(0,0){\circle*{4}} \put(0,40){\circle*{4}}
 \put(0,40){\line(0,-1){40}}
 \put(40,0){\circle*{4}} \put(40,40){\circle*{4}}
 \put(80,0){\circle*{4}} \put(80,40){\circle*{4}}
 \put(120,0){\circle*{4}} \put(120,40){\circle*{4}}
 \put(80,40){\line(0,-1){40}}
 \put(0,0){\line(1,0){120}} \put(0,40){\line(1,0){120}}
 \put(120,40){\line(0,-1){40}}
 \put(150,0){\circle*{4}} \put(150,40){\circle*{4}}
 \put(150,40){\line(0,-1){40}}
 \put(150,0){\line(1,0){40}} \put(150,40){\line(1,0){40}}
 \put(190,0){\circle*{4}} \put(190,40){\circle*{4}}
 \put(230,20){\circle*{4}}
 \put(230,20){\line(-2,-1){40}} \put(230,20){\line(-2,1){40}}
 \put(230,20){\line(2,-1){40}} \put(230,20){\line(2,1){40}}
 \put(270,0){\circle*{4}} \put(270,40){\circle*{4}}
 \put(270,40){\line(0,-1){40}}
 \put(230,31){\circle{20}}
\end{picture}
\end{center}
\vskip -0.2cm \caption{Edge-pinching a chordal even cycle}
\label{figure_chordal_evengon}
\end{figure}

\begin{Corollary}\label{polarizability_not_invariant}
Polarizability is not an invariant property of the algebra $A$.
\end{Corollary}

\demo Consider $\ff$ such that the graph ${\cal G}(\ff)$ is an even
cycle with one single chord inducing a decomposition in two smaller
even cycles and $\ff'$ whose associated graph ${\cal G}(\ff')$ is
obtained by edge-pinching the chord of ${\cal G}(\ff)$. Then $\ff$
is polarizable by Theorem~\ref{polar_characterization} while $\ff'$
is not polarizable by Proposition~\ref{pinching_chordal_evengon}.
Nevertheless, $k[\ff]\simeq k[\ff']$ (the defining ideals of both
$k$-subalgebras coincide).
 \qed

 \begin{Remark}\label{polarizability_and_isomorphism}\rm
 The actual reason why polarizability is not an invariant property of
 the algebra $A$ is that a $k$-algebra isomorphism may not preserve
 certain crucial configurations of the corresponding graph.
 Thus, e.g., in Figure~\ref{figure_chordal_evengon} the
 path-degenerate bow tie in the right most graph, whose odd cycles
 are a pentagon and a triangle, is not preserved under the above
 isomorphism of algebras.
 \end{Remark}

\subsection{Polarizability versus normality}

Recall the notion of a cohesive set of monomials.

\begin{Definition}[{\cite[Definition 4.2]{SimisVilla}}]\label{defining_cohesive}\rm
The set $\ff$ is said to be {\it cohesive} if there is no partition
$\xx=\yy\cup\zz$ of the variables such that $\ff=\bg\cup\hh$, where
the monomials in the set $\bg$, resp. $\hh$, involve only the
$\yy$-variables, resp. $\zz$-variables.
\end{Definition}

One clearly has that $\ff$ is cohesive if and only if ${\cal
G}(\ff)$ is connected. The following characterization for the
normality of $k[\ff]$ has essentially been obtained (independently)
in \cite{bowtie} and \cite{ohsugihibi}.

\begin{Proposition}\label{characterization_normality}
Let $A=k[\ff]\subset R$ be generated by a cohesive set $\ff$ of
monomials of degree $2$ and let ${\cal G}(\ff)$ denote the
corresponding graph. The following conditions are equivalent:
\begin{enumerate}
\item[{\rm (1)}]
$A$ is integrally closed;
\item[{\rm (2)}] ${\cal G}(\ff)$ satisfies the so-called odd cycle
condition, i.e.,  for any two odd cycles which are induced {\rm
(}i.e., no chords{\rm)} in ${\cal G}(\ff)$ and have mutually
disjoint vertex sets, there exists an edge of ${\cal G}(\ff)$
joining a vertex of one cycle to a vertex of the other.
\item[{\rm (3)}]
Any induced bow tie of ${\cal G}(\ff)$ is either a path-degenerate
bow tie or a monedge bow tie {\rm (}possibly including the
respective looped versions{\rm )}.
\end{enumerate}
\end{Proposition}

\demo (1) $\Leftrightarrow$ (2) is (i) $\Leftrightarrow$ (iii) in
\cite[Corollary~2.3]{ohsugihibi}.

\noindent (2) $\Rightarrow$ (3) This is obvious.

\noindent (3) $\Rightarrow$ (2) Given two odd cycles as stated --
called for convenience {\it non-chordal\/} -- there must be a path
connecting the two since we are assuming that ${\cal G}(\ff)$ is
connected. This yields a bow tie ${\cal B}$ in the graph, and we may
assume that ${\cal B}$ has a connecting path of smallest length
$\ell$ among all bow ties in the graph whose structural odd cycles
are non-chordal. Assume, as if it were, that $\ell\geq 2$. If ${\cal
B}$ is induced, it would be a contradiction to (3). If it is not
induced, let $e$ be an edge between two vertices of ${\cal B}$.
Since the two odd cycles are non-chordal, $e$ must connect vertices
across the two cycles or across a cycle and the path. In the first
case, we are done, while the second case is ruled out as it implies
a new bow tie with non-chordal cycles such that $e$ is an edge of
one of the cycles and admitting a connecting path of length $\leq
\ell-1$. \qed

\medskip

The next result explains the precise relationship between the
notions of polarizability and normality.

\begin{Theorem}\label{polar_is_normal}
Let $A=k[\ff]\subset R$ be generated by a cohesive set $\ff$ of
monomials of degree $2$ and let ${\cal G}(\ff)$ denote the
corresponding graph.
\begin{enumerate}
\item[{\rm (i)}] If $\ff$ is polarizable then $A$ is integrally closed {\rm
(}hence, a Cohen--Macaulay ring{\rm )}.
\item[{\rm (ii)}] Conversely, suppose that ${\cal G}(\ff)$ has no configuration of the
following kinds:
\begin{itemize}
\item[{\rm (a)}] Induced monedge bow ties {\rm(}with neither odd cycle
degenerating into a loop{\rm )}$\,${\rm ;}
\item[{\rm (b)}] Induced monedge looped bow ties {\rm(}with only one odd cycle
degenerating into a loop{\rm )}$\,${\rm ;}
\item[{\rm (c)}] Induced path-degenerate bow ties {\rm(}with neither odd cycle
degenerating into a loop{\rm )}.
\end{itemize}
If $A$ is integrally closed then $\ff$ is polarizable.
\end{enumerate}
\end{Theorem}
\demo It follows from Proposition~\ref{characterization_normality}
and Theorem~\ref{polar_characterization}. \qed

\begin{Remark}\rm
Note that the above result does not conflict with the result of
Corollary~\ref{polarizability_not_invariant} (see also
Remark~\ref{polarizability_and_isomorphism}).
\end{Remark}

The following consequence for algebras of Veronese type of degree
$2$ could have been given before with slightly more effort, but
having it here stresses the normality of these algebras. Recall
that, given an integer $d\geq 1$ and a sequence of integers $1\leq
s_1\leq\cdots\leq s_n\leq d$, the $k$-subalgebra $A\subset R$
generated by the set of monomials
$$F=\{x^{a_1}\cdots x^{a_n}\,|\,
a_1+\cdots+a_n=d;\
 0\leq a_i\leq s_i\, \forall\, i\}$$
is called the {\it algebra of Veronese type\/} of degree $d$
subordinate to the vector $(s_1,\ldots,s_n)$. These algebras form a
subclass of the class of the polimatroidal algebras of maximal rank
(see \cite{DeHi}, \cite{HeHi}). In the next subsection we will
actually show that {\it all\/} polimatroidal algebras of degree $2$
are polarizable.

\begin{Corollary}\label{veronese_type}
If $\ff$ are the defining generators of an algebra $A\subset R$ of
Veronese type of degree $2$ then $\ff$ is polarizable.
\end{Corollary}
\demo It is known that $A$ is normal (cf., e.g., \cite{DeHi}; see
also \cite{bowtie}). On the other hand, since $d=2$ the relevant
subordinating vectors have $s_i\leq 2$ for all $i$. It follows that
$\ff$ consists of all squarefree monomials of degree $2$ and
possibly some pure powers. It is then self-evident that the
associated graph does not admit any induced path-degenerate or
monedge bow ties except eventually looped-triangles or two loops
joined with an edge. By Theorem~\ref{polar_is_normal}, (ii), $\ff$
is polarizable. \qed

\begin{Corollary}\label{polar_is_birational}
Let $F:\pp^{n-1}\dasharrow \pp^{m-1}$ be a rational map defined by a
cohesive set $\ff$ of distinct monomials of degree $2$. If $\,\dim
k[\ff]=n$ and $\ff$ is polarizable then $F$ maps $\pp^{n-1}$
birationally onto its image. In particular, $k[\ff]$ is a rational
singularity.
\end{Corollary}
\demo First observe that the claim on birationality is equivalent to
saying that the ring extension $k[\ff]\subset k[(\xx)_2]$
($2$-Veronese) is birational (see, e.g., \cite[Proof of Proposition
2.1]{SimisVilla0}). Thus, if for some subset $\ff'\subset \ff$ the
corresponding rational map is birational onto its image then so will
be the one defined by $\ff$. Let us choose $\ff'$ to be the subset
of the squarefree monomials in $\ff$.

Now, on one hand Theorem~\ref{polar_is_normal}, (i), implies that
$k[\ff']$ is normal, while on the other hand, the normality of the
squarefree  $k[\ff']$ is equivalent to the normality of the ideal
$(\ff')$ in this case (see \cite[Corollary 8.7.13]{VillaBook}).
Therefore, by \cite[Proposition 3.1]{SimisVilla0} the extension
$k[\ff']\subset k[(\xx)_2]$ is birational, as required. \qed

\subsection{Polarizability versus linear presentation}

We deal here with the case in which $\ff$ is linearly presented,
i.e. when its module of first syzygies is generated by linear ones.
We characterize this property in terms of the diameter of a graph
(Lemma~\ref{lin_pres_char}) and show that if $\ff$ is linearly
presented then it is polarizable
(Proposition~\ref{lin_pres_is_polarizable}).

\medskip

In order to characterize when $\ff$ is linearly presented, we
introduce the {\it edge graph} of ${\cal G}(\ff)$, denoted ${\cal
L}(\ff)$ (see \cite[Definition~6.6.1]{VillaBook}): its vertex set is
the set of edges of ${\cal G}(\ff)$, hence can be viewed as the
elements of $\ff$; two vertices $f_i$ and $f_j$ of ${\cal L}(\ff)$
are adjacent (i.e., form an edge) if and only if $f_i$ and $f_j$
have a common variable (i.e., $\gcd{(f_i,f_j)}\neq 1$). Observe that
the graph ${\cal L}(\ff)$ is always a simple graph (no loops) and
that $\ff$ is cohesive (see Definition~\ref{defining_cohesive}
previously recalled) if and only if ${\cal G}(\ff)$ is connected, if
and only if ${\cal L}(\ff)$ is connected.

\begin{Example}\rm
For $\ff=\{f_1,\ldots,f_4\}$ with $f_1=x_1^2$, $f_2=x_1x_2$,
$f_3=x_2x_3$ and $f_4=x_1x_3$, the graphs ${\cal G}(\ff)$ and ${\cal
L}(\ff)$ are given in Figure~\ref{figure_graph_edgegraph}.

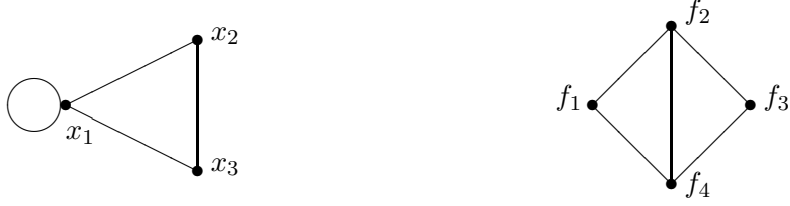
\begin{figure}[h]
\begin{center}
\setlength{\unitlength}{.035cm}
\begin{picture}(0,50)(0,10)
\put(-130,30){\circle*{4}} \put(-130,17){$x_1$}
\put(-142,30){\circle{20}} \put(-130,30){\line(2,1){50}}
\put(-130,30){\line(2,-1){50}} \put(-80,55){\line(0,-1){50}}
\put(-80,55){\circle*{4}} \put(-75,55){$x_2$}
\put(-80,5){\circle*{4}} \put(-75,5){$x_3$}
\put(100,60){\line(0,-1){60}} \put(100,60){\line(1,-1){30}}
\put(100,0){\line(1,1){30}} \put(100,60){\line(-1,-1){30}}
\put(100,0){\line(-1,1){30}} \put(100,0){\circle*{4}}
\put(105,-3){$f_4$} \put(100,60){\circle*{4}} \put(105,63){$f_2$}
\put(130,30){\circle*{4}} \put(135,30){$f_3$}
\put(70,30){\circle*{4}} \put(56,30){$f_1$}
\end{picture}
\end{center}
\vskip -0.2cm \caption{A graph and its edge
graph}\label{figure_graph_edgegraph}
\end{figure}
\end{Example}

As observed in \cite[Lemma~4.1]{SimisVilla}, the lack of
cohesiveness is an obstruction for the existence of enough linear
syzygies. In the situation we focus on in this section, it is thus
natural to assume that $\ff$ is cohesive, i.e., that ${\cal G}(\ff)$
and ${\cal L}(\ff)$ are both connected graphs.

\begin{Definition}\rm
Given a simple connected graph ${\cal G}$, the {\it distance\/}
between two vertices of ${\cal G}$ is the minimum length of a path
connecting them, and the {\it diameter\/} of ${\cal G}$ is the
longest distance (i.e., the longest shortest path) between any two
of its vertices.
\end{Definition}

\begin{Lemma}\label{lin_pres_char}
Assume that $\ff$ is cohesive. Then, the ideal $I=(\ff)\subset R$ is
linearly presented if and only if the graph ${\cal L}(\ff)$ is of
diameter $\leq 2$.
\end{Lemma}
\demo Recall that ${\bf f}=\{f_1,\ldots,f_m\}$, denote by
$\{e_1,\ldots,e_m\}$ the canonical basis of the free module $R^m$,
and  set
$$s_{ij}:=\frac{f_j}{\gcd(f_i,f_j)}\,e_i-
\frac{f_i}{\gcd(f_i,f_j)}\,e_j\in R^m,$$
 for $i,j\in\{1,\ldots,m\}$,
$i\neq j$. It is well-known (see, e.g., \cite[Chapter~5,
Thm.~3.2]{CoxUsingAG}) that the first sygygy module of $I$ is
generated by the set ${\mathcal S}(\ff):=\{s_{ij}\,|\,1\leq i<j\leq
m\}$. Consider the partition ${\mathcal
S}(\ff)=\mathcal{LS}(\ff)\cup \mathcal{KS}(\ff)$ where
$$\mathcal{LS}(\ff):=\{s_{ij}\,|\, 1\leq i<j\leq m,\;
\gcd(f_i,f_j)\neq 1\}$$
and
$$ \mathcal{KS}(\ff):=\{s_{ij}\,|\,1\leq
i<j\leq m,\; \gcd(f_i,f_j)=1\}.$$
 The syzygies $s_{ij}$ in $\mathcal{LS}(\ff)$ are linear, and the ones in
 $\mathcal{KS}(\ff)$
 are Koszul syzygies since $s_{ij}=f_je_i-f_ie_j$ if $\gcd(f_i,f_j)=1$. The
ideal $I=(\ff)$ has linear syzygies if and only if
$\mathcal{KS}(\ff)$ is contained in the submodule of $R^m$ generated
by $\mathcal{LS}(\ff)$.

\smallskip
First observe that the diameter of the graph ${\cal L}(\ff)$ is 1
(i.e., the graph ${\cal L}(\ff)$ is complete) if and only if
$\mathcal{KS}(\ff)=\emptyset$. More precisely, for all $i,j$, $1\leq
i<j\leq m$, one has that the distance between the vertices $f_i$ and
$f_j$ of ${\cal L}(\ff)$ is $1$ if and only if
$s_{ij}\in\mathcal{LS}(\ff)$.

The result will follow if one shows that if
$\mathcal{KS}(\ff)\neq\emptyset$ then, for any $i,j$ such that
$s_{ij}\in\mathcal{KS}(\ff)$, the syzygy $s_{ij}$ belongs to the
submodule generated by $\mathcal{LS}(\ff)$ if and only if the
distance between the vertices $f_i$ and $f_j$ of ${\cal L}(\ff)$ is
$2$.

\smallskip
Thus, suppose $\mathcal{KS}(\ff)\neq\emptyset$ and let
$g\in\mathcal{KS}(\ff)$. One can assume, without loss of generality,
that $g=s_{12}$ and, relabelling the variables if necessary, that
$f_1=x_1x_i$ and $f_2=x_jx_n$ for some $i\in\{1,\ldots,n-1\}$ and
$j\in\{2,\ldots,n\}$ such that $i\neq j$. Then,
$g=x_jx_ne_1-x_1x_ie_2$. If $g$ belongs to the submodule generated
by $\mathcal{LS}(\ff)$, then there exists at least one element in
$\mathcal{LS}(\ff)$ such that one of its two nonzero entries is
either $x_je_1$ or $x_ne_1$. This implies that either
$x_jx_1\in\ff$, or $x_jx_i\in\ff$, or $x_nx_1\in\ff$, or
$x_nx_i\in\ff$, and hence, the distance between the vertices $f_1$
and $f_2$ of ${\cal L}(\ff)$ is 2. Conversely, if the distance
between the vertices $f_1$ and $f_2$ of ${\cal L}(\ff)$ is 2, one
has that either $x_jx_1\in\ff$, or $x_jx_i\in\ff$, or
$x_nx_1\in\ff$, or $x_nx_i\in\ff$. Assume for example that
$f_3=x_jx_1$. Then, $s_{13}=x_je_1-x_ie_3$ and
$s_{23}=x_1e_2-x_ne_3$ are elements in $\mathcal{LS}(\ff)$, and
since $g=x_ns_{13}-x_is_{23}$, we are through. \qed

\begin{Remark}\label{froberg_char_lin_resol}\rm
There is another kind of complementary configuration to a given
simple graph ${\cal G}(\ff)$ called the {\em complement} of ${\cal
G}(\ff)$, denoted $\overline{{\cal G}(\ff)}$: it has the same vertex
set as ${\cal G}(\ff)$, and the edges are those edges of the
complete simple graph on the same vertex set which are not edges of
${\cal G}(\ff)$ (see \cite[~p.~175]{VillaBook}).

Fr\"oberg (\cite{Froberg}) proved that the ideal $I=(\ff)\subset R$
generated by a set $\ff$ of square-free monomials of degree $2$ has
a linear resolution if and only if the graph $\overline{{\cal
G}(\ff)}$ is chordal, i.e., has no induced cycles of length $\geq
4$. This result is related to Lemma~\ref{lin_pres_char} in the
following way:  if $\ff$ is a set of square-free monomials of degree
$2$, the graph ${\cal L}(\ff)$ has diameter $\leq 2$ if and only if
the graph $\overline{{\cal G}(\ff)}$ has no induced $4$-cycles.
Thus, for simple graphs Lemma~\ref{lin_pres_char} reproves a piece
of Fr\"oberg's result. Actually, there is a refinement of
Fr\"oberg's result in \cite[Theorem~2.1]{Eisetal} which we regrettably have been unaware of.
Using it together with
\cite[Proposition~2.3]{Eisetal}, one can recover
Lemma~\ref{lin_pres_char}. Since the above proof is straightforward
and elementary, we decided to keep it (see also \cite{GimOsc} for
yet another approach).
\end{Remark}

We can now prove the following fundamental connection between linear presentation
and polarizability.

\begin{Proposition}\label{lin_pres_is_polarizable}
If the ideal $I=(\ff)\subset R$ generated by a set $\ff$ of
monomials of degree $2$ is linearly presented then $\ff$ is
polarizable.
\end{Proposition}
\demo By the characterization in Lemma~\ref{lin_pres_char}, if
$I=(\ff)\subset R$ is linearly presented, the induced odd cycles
(with no chord) in ${\cal G}(\ff)$ (if any) are loops and triangles.
Moreover, the induced bow ties in ${\cal G}(\ff)$ (if any) are two
loops joined with an edge or a triangle with a loop centered in one
of its vertices. By Theorem~\ref{polar_characterization}, $\ff$ is
polarizable. \qed

\begin{Corollary}\label{polimatroidal_is_polar}
If $\ff$ is a polimatroidal set of monomials of degree $2$ then
$\ff$ is polarizable.
\end{Corollary}
\demo By \cite{ConHer}, if $\ff$ is ordered in the reverse
lexicographic order, then it has linear quotients, i.e., the ideals
$(f_1,\ldots,f_{i-1}):f_i$ are generated by a set of variables, for
every $1\leq i\leq m$. It is self-evident that having linear
quotients entice linear presentation, hence the result follows from
Proposition~\ref{lin_pres_is_polarizable}. \qed

\medskip

In a curious roundabout fashion we recover \cite[Corollary
3.8]{SimisVilla}:

\begin{Corollary}\label{linear_pres_is_birational}
Let $F:\pp^{n-1}\dasharrow \pp^{m-1}$ be a rational map defined by a
cohesive set $\ff$ of distinct monomials of degree $2$. If $\dim
k[\ff]=n$ and $(\ff)\subset R$ is linearly presented then $F$ maps
$\pp^{n-1}$ birationally onto its image. In particular, $k[\ff]$ is
a rational singularity.
\end{Corollary}
\demo It follows immediately from
Proposition~\ref{lin_pres_is_polarizable} and
Corollary~\ref{polar_is_birational}. \qed

\bigskip

We end with a couple of remarks and an example.

Namely, let again $k[\TT]/P\simeq k[\ff]$. If $P$ happens to be
generated by sole quadrics then a minimal set of generators of the
polar syzygy module $\fp$ is automatically a minimal subset of
generators of the differential syzygy module $\fz$. This is of
course a favorable situation which one would like to understand
better. If the ideal $(\ff)\subset R$ is linearly presented then $P$
has ``many'' quadrics, but still may require generators of higher
degrees. In fact, these degrees may be arbitrarily high as the
following example shows.

\begin{Example}\label{bounding_cycle}\rm Consider a complete graph (no loops) with $t\geq 3$
vertices. Mark one of the $t$-cycles of the graph as the ``bounding
cycle''. For each pair of consecutive vertices $v_1,v_2$ of the
bounding cycle introduce a new vertex $v$ and new edges $vv_1$ and
$vv_2$. In this way we have constructed a graph on $n=2t$ vertices
equipped with a new bounding $n$-cycle. It is easy to see that the
diameter of the new graph is still $\leq 2$, hence its edges
correspond to a set $\ff$ that is linearly presented, hence
polarizable by Proposition~\ref{lin_pres_is_polarizable}. However
the new bounding cycle induces an element of $P$ of degree $t$ that
is not contained in the ideal generated by the quadrics in $P$. The
reason it does not induce an extra minimal generator at the level of
$\fz$ (or, which is the same, of $\fp$) is that it is decomposed by
the internal chords of the new bounding cycle.
\end{Example}

Next is depicted the simplest case ($t=3$).

\vspace{2cm}

\begin{figure}[h]
\setlength{\unitlength}{.035cm}
\begin{picture}(250,50)
\multiput(150,70)(40,20){2}{\circle*{4}} \put(230,70){\circle*{4}}
\multiput(150,20)(80,0){2}{\circle*{4}} \put(190,0){\circle*{4}}
\put(150,70){\line(2,1){40}} \put(230,70){\line(-2,1){40}}
\put(150,20){\line(0,1){50}} \put(230,20){\line(0,1){50}}
\put(190,0){\line(2,1){40}} \put(190,0){\line(-2,1){40}}
\put(150,70){\line(2,0){80}} \put(190,0){\line(3,5){40}}
\put(150,70){\line(3,-5){40}} \put(190,100){$1$} \put(240,70){$2$}
\put(240,10){$3$} \put(190,-15){$4$} \put(140,10){$5$}
\put(140,70){$6$}
\end{picture}
\vskip 0.5cm \caption{Linearly presented with cubic relations}
\end{figure}
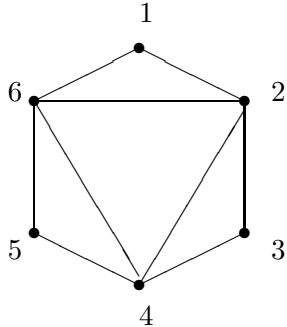

A question also naturally arises as to what is the impact on
polarizability of $\ff$ if the presentation ideal $P$ is actually
fully generated in degree $2$. Easy examples show that, in general,
$\ff$ may not be polarizable. However, these examples are such that
the ideal $P\subset k[\TT]$ is not itself linearly presented. Thus
it seems reasonable to pose:

\begin{Question}\rm
Suppose that $P$ is generated by quadrics and is linearly presented.
Is $\ff$ polarizable? More strongly, is $\ff$ linearly presented as
well?
\end {Question}

A special important class of algebras satisfying these hypotheses
are the Koszul algebras $A=k[\TT]/P$, which are generated by
quadrics and have linear resolution.

{\small

\bigskip

{\sc Isabel Bermejo}, Facultad de Matem\'aticas, Universidad de La
Laguna, 38200 La Laguna, Tenerife, Canary Islands, Spain

{\it Email}: {\tt ibermejo@ull.es}

\bigskip

{\sc Philippe Gimenez}, Departamento de Algebra, Geometr\'{\i}a y
Topolog\'{\i}a, Facultad de Ciencias, Universidad de Valladolid,
47005 Valladolid, Spain

{\it Email}: {\tt pgimenez@agt.uva.es}

\bigskip

{\sc Aron Simis}, Departamento de Matem\'atica, CCEN, Universidade
Federal de Pernambuco, Cidade Universit\'aria, 50740-540 Recife, PE,
Brazil

{\it Email}: {\tt aron@dmat.ufpe.br}

}


\begin{thebibliography}{99}



\bibitem{ConHer}{A. Conca and J. Herzog, Castelnuovo-Mumford regularity of
products of ideals, Collect. Math. {\bf 54} (2003), 137-152.}

\bibitem{CoxUsingAG}{D. Cox, J. Little and D. O'Shea, {\it Using Algebraic Geometry}, Springer,
New York-Berlin-Heidelberg, 1998.}

\bibitem{Eisetal}{D. Eisenbud, M. Green, K. Hulek and S. Popescu,
Restricting linear syzygies: algebra and geometry, Compos. Math.
{\bf 141} (2005), 1460-1478.}

\bibitem{GimOsc}{O. Fern\'andez-Ramos and P. Gimenez, Nonlinear syzygies of smallest degree of an ideal
associated to a graph, preprint 2008.}

\bibitem{Froberg}{R. Fr\"oberg, On Stanley-Reisner rings. In: {\it Topics in
Algebra, Part~2} (Warsaw, 1988), Banach Center Publ. {\bf 26} (1990)
57--70.}

\bibitem{GN}{P. Gordan und M. Noether, Ueber die algebraischen Formen, deren
Hesse'sche Determinante identisch verschwindet, Math. Ann., {\bf 10}
(1876), 547--568.}

\bibitem{HeHi}{J. Herzog and T. Hibi, Discrete polymatroids,
J. Algebraic Combin. {\bf 16} (2002), 239–-268.}

\bibitem{DeHi}{E. de Negri and T. Hibi, Gorenstein algebras of Veronese type,
J. Algebra {\bf 193} (1997) 629--639.}

\bibitem{ohsugihibi}{H. Ohsugi and T. Hibi, Normal polytopes arising from finite graphs,
J. Algebra {\bf 207} (1998) 409--426.}

\bibitem{jac}{A. Simis, On the Jacobian module associated to  graph,
Proc. Amer. Math. Soc. {\bf 126} (1998), 989--997.}

\bibitem{SimisDiff}{A. Simis, Two differential themes in characteristic zero,
{\it in\/} {\sc Topics in Algebraic and Noncommutative Geometry},
Proceedings in Memory of Ruth Michler (Eds. C. Melles, J.-P.
Brasselet, G. Kennedy, K. Lauter and L. McEwan), {\sc Contemporary
Mathematics} {\bf 324}, Amer. Math. Soc., Providence, RI, 2003,
195--204.}

\bibitem{itg}{A. Simis, W. V. Vasconcelos and R. Villarreal, On the ideal theory of graphs,
J. Algebra {\bf 167} (1994), 389--416.}

\bibitem{bowtie}{A. Simis, W. V. Vasconcelos and R. Villarreal, The integral closure
of subrings associated to graphs, J. Algebra {\bf 199} (1998),
281--289.}

\bibitem{SimisVilla0}{A. Simis and R. Villarreal, Constraints for the normality of
monomial subrings and birationality, Proc. Amer. Math. Soc. {\bf
131} (2003), 2043-2048.}

\bibitem{SimisVilla}{A. Simis and R. Villarreal, Linear syzygies and birational
combinatorics, Results Math. {\bf 48} (2005), 326--343.}

\bibitem{VillaBook}{R. H. Villarreal, {\it Monomial
Algebras\/}, Monographs and Textbooks in Pure and Applied Mathematics {\bf 238},
Marcel Dekker, New York, 2001.}

\end{thebibliography}
\end{document}